\documentclass[a4paper,10pt]{scrartcl}

\usepackage{amsmath}
\usepackage[amsthm]{ntheorem}
\usepackage{amssymb}
\usepackage{color} 
\usepackage{epsfig}
\usepackage{graphicx}
\usepackage{subfigure}
\usepackage{booktabs}
\usepackage{lscape}
\usepackage{dcolumn}
\usepackage{supertabular}
\usepackage{pstricks}
\usepackage{pst-node}
 \usepackage{pst-tree}
\usepackage{bbm}
\usepackage[latin1]{inputenc}
\usepackage{psfrag}

 \usepackage{hyperref} 
 \usepackage{algorithm}

	\theoremstyle{break}
	\newtheorem{satz}{Theorem}[section]

	\newtheorem{bem}[satz]{Remark}
	
\newcommand{\dis}{\displaystyle}

\def\ZZ{{\mathbb Z} }
\def\RR{{\mathbb R} }

\def\EE{{\mathbb E} }
\def\dis{\displaystyle}

\def\Cor{{\mathbb Cor} }
\def\Cov{{\mathbb Cov} }
\def\Id{\text{Id}}

\title{Multi-Level Monte Carlo approaches for numerical homogenization}
\author{Yalchin Efendiev\thanks{Department of Mathematics, Texas A \& M
    University, College Station, TX 77845, USA} \and Cornelia
  Kronsbein\thanks{Fraunhofer ITWM \& University of Kaiserslautern,
    67663 Kaiserslautern, Germany} \and Fr\'ed\'eric
  Legoll\thanks{Laboratoire Navier, \'Ecole des Ponts ParisTech,
Universit\'e Paris-Est, 6 et 8 avenue Blaise Pascal, 77455
Marne-La-Vall\'ee Cedex 2, France and 
INRIA Rocquencourt, MICMAC Team-Project, Domaine de Voluceau, B.P. 105,
78153 Le Chesnay Cedex, France}}
\date{\today}
\newcolumntype{d}{D{.}{.}{7}}

\begin{document}
\maketitle
\abstract{{\bf Abstract.}
In this article, we study the application of Multi-Level Monte
Carlo (MLMC) approaches to numerical random homogenization. Our objective
is to compute 
the expectation of some functionals of the homogenized coefficients,
or of the homogenized solutions. This is accomplished within MLMC by
considering different levels of representative volumes (RVE), and, when
it comes to homogenized solutions, different levels of coarse-grid
meshes. Many inexpensive computations with the smallest
RVE size  and the largest coarse mesh are combined with fewer expensive
computations performed on larger RVEs and smaller coarse meshes. We show
that, by carefully selecting the number of realizations at each level,
we can achieve a speed-up in the computations in comparison to a
standard Monte Carlo method. 
Numerical results are presented both for one-dimensional
and two-dimensional test-cases.
}

\section{Introduction}

Many multi-scale problems have uncertainties at the smallest scales, that
are due to the incomplete knowledge one has of the
microstructure. 
For example, when considering porous materials, the microstructure
is often generated based on some limited statistical information.
This can lead to large uncertainties in terms of microscale heterogeneities.
These uncertainties at micro-scales need to be mapped onto the simulations
on a coarse-grid, and this typically leads to considering large
representative volumes (RVE) for these microstructures.

In practice, the upscaled quantities that are used at the
macroscopic level are computed using the solution of some local problems
posed on these microstructures. It is often needed to solve many such
local problems (corresponding to many different random realizations, or
snapshots, of the microstructure), each of which being expensive due to the presence of small
scales. The resulting amount of computational work may thus be
prohibitively expensive. In this article, our objective 
is to design a computational approach that allows for fast calculations
of the coarse-scale quantities based on fewer realizations. 

Our idea is to apply the Multi-Level Monte Carlo (MLMC) framework to
multi-scale simulations. The MLMC approach was first introduced by
Heinrich in~\cite{H01} for finite- and infinite-dimensional integration. 
Later on, it was applied to stochastic ODEs by Giles
(see~\cite{G08,G08b}). More recently, this approach has been used for
PDEs with stochastic coefficients by several authors,
see~\cite{BSZ10,CGST11,abdulle,charrier,teckentrup}.
To compute an approximation of the expectation $\EE(X)$ of some random
variable $X$, the MLMC approach consists in considering several random
variables $X_l$, at different levels $l$, that approximate $X$ with
various accuracies.
The main idea is then to use different numbers of samples
(i.e. independent realizations) at different levels.
More precisely, many samples are used at the coarsest, less accurate level
where the computation for each realization is
inexpensive, while fewer samples are
used at the finest, most accurate level that is expensive to compute.
Combining the results of these computations by
carefully selecting the number of realizations at each level
can speed-up the computations in comparison to a standard Monte Carlo
(MC) approach, where only one level (that of the quantity of interest
itself) is considered. See Section~\ref{sec:mlmc_gene} below for more
details on the MLMC approach.

\medskip

In the framework of numerical stochastic homogenization, local problems
are solved on representative volumes (RVE), and apparent effective properties
are next defined as averages of the solutions of these local problems
over the RVEs.
The computations on the RVEs are usually expensive, because large RVEs
need to be considered to obtain effective properties with a reasonable
accuracy. In the framework of MLMC approaches, our idea is to use RVEs
of different sizes, and to consider many independent realizations of the
smaller ones, for which the associated local problem is inexpensive to
solve, and fewer realizations of the larger ones.

The convergence of the MLMC approach depends on the
accuracy of the computations at each level. Assessing how this
accuracy improves when more expensive computations are considered
is critical to determine how to choose the number of realizations at
each level. In our case, we thus have to determine how the accuracy of
apparent effective properties depend on the RVE size. Such estimations
are not easy to obtain, both from a theoretical and a practical viewpoint. 
In this work, we use the fact that, under some assumptions on the heterogeneous
coefficients, it is known that the accuracy of the effective property
approximation scales as $(\epsilon/\eta)^\beta$ for some $\beta>0$, where
$\eta$ is the RVE size and $\epsilon$ is the characteristic small
lengthscale of the heterogeneities  
(see
e.g.~\cite{B07,BCLL11_mprf,BCLL11_lncse,BP04,CLL10,E99diss,GO10,JKO94,Y87}).

When the MLMC approach is used to compute the expectation of some
functionals of the homogenized {\em solution} (rather than the
homogenized {\em coefficient}), we can use RVEs of different size to compute
the homogenized coefficients, and also coarse grids with various size to
solve the coarse scale equation. In addition to assessing the accuracy
of the approximation of the effective properties in each RVE, we need to
assess the accuracy when solving the
coarse-scale equation. Standard FEM results are then useful. 

An important remark is that MLMC approaches are interesting when
effective properties are stochastic (otherwise, such approaches are as
efficient as a standard MC approach). This situation appears
in many applications, although homogenization theories
for this case are less studied. Most homogenization theories
are indeed developed for ergodic coefficients that vary over a single
scale. In this case, the apparent homogenized quantities, when computed
on infinitely large RVEs, are deterministic. In the sequel, we briefly
discuss homogenization results when the homogenized coefficient is
stochastic (even when infinitely large RVEs are considered), and we use
these results in our MLMC approach to adequately select the number of
realizations at each level (namely, for each RVE size and each coarse
grid size).

\medskip

Consider now the specific question of computing homogenized solutions with
several grids of different size. For each of these grids, we first need
to precompute the effective properties, say at each Gauss point of the
macroscopic grid. Assume that these coarse grids are nested. Then, once the
effective properties have been computed at the finest level (i.e. for the
Gauss points of the finest grid), no additional precomputation is needed
to compute effective properties for the coarser grids (since their
Gauss points are a subset of the Gauss points of the finest grid). 
In this case, we propose to use a weighted MLMC approach, where we give
different weights to each level, so as to optimize the accuracy at a
given cost.

\medskip

Our article is organized as follows. In Section~\ref{sec:prelim}, we
briefly review theoretical homogenization results and describe in
details the MLMC approach in a general context. In
Section~\ref{sec:mlmc}, we next describe how to apply the MLMC approach
to compute an approximation of the homogenized coefficients, and assess
the accuracy of the proposed approach. 
We next turn in Section~\ref{sec:homog_sol} to the
computation of the homogenized solutions, using either the MLMC
or the weighted MLMC approaches.
Numerical results are collected in Section~\ref{sec:num}. We consider
the case of the effective coefficients in
Section~\ref{sec:num_coeff}, and of the homogenized solutions in
Section~\ref{sec:num_sol}. In both cases, we show that the MLMC approach
yields a significant speed-up in comparison to a standard MC approach.

\section{Preliminaries}
\label{sec:prelim}

\subsection{Numerical homogenization}
\label{numHom}

In this section, we describe the numerical homogenization procedure we
use. Consider the problem
\begin{equation}
\label{microprob}
-\text{div}(A_\epsilon(x,\omega)\nabla u_\epsilon)=f \text{ in $D$},
\end{equation}
where $D$ is an open bounded subset of $\RR^d$, 
$A_\epsilon(x,\omega)$ is a heterogeneous random field (with a small
characteristic length scale $\epsilon$), $\omega$
designates a random realization and $f\in L^2(D)$ is a non-random
function. We complement the problem~\eqref{microprob} with some boundary
conditions that we do not specify, such that its solution
$u_\epsilon$ is well defined (for instance, $u_\epsilon = 0$ on
$\partial D$ almost surely). Furthermore, we assume that $A_\epsilon$ is
uniformly bounded and coercive, in the sense that there exists two
positive deterministic numbers $0 < a_{\text{min}} \leq a_{\text{max}}$
such that, for any $\epsilon$, any $\xi \in \RR^d$ and any $1 \leq i,j
\leq d$, 
$$
a_{\text{min}} | \xi |^2 \leq 
\xi^T A_\epsilon(x,\omega) \xi,
\quad
\left| \left[ A_\epsilon(x,\omega)
\right]_{ij} \right| \leq a_{\text{max}},
$$
almost everywhere in $D$ and almost surely.

For almost all realizations $\omega$, we consider a numerical
homogenization procedure as follows. Given a representative volume
centered at a macroscopic point $x$ with size $\eta$, 
$$
Y_{\eta}^x = \left( x- \frac{\eta}{2},x+\frac{\eta}{2} \right)^d,
$$ 
we solve, for any $1 \leq i \leq d$, the local problems 
\begin{equation}
\label{eq:corrector}
\text{div}(A_\epsilon(y,\omega)\nabla \chi_i(y,\omega)) = 0 
\text{ in $Y_{\eta}^x$},
\quad
\chi_i(y,\omega) = y_i \ \text{on $\partial Y_{\eta}^x$}.
\end{equation}
Note that the precise boundary conditions used in these local problems
are not essential when there 
is a scale separation. Rather than Dirichlet boundary conditions as
in~\eqref{eq:corrector},
it is also possible to use Neumann boundary conditions, or periodic
boundary conditions (see~\cite{BP04,Kanit2003}).

\begin{figure}[tbp]
\centering
\scalebox{0.4}{\input{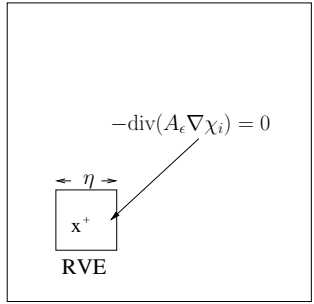}}
\caption{Illustration of the numerical homogenization procedure.}
\label{rev}
\end{figure}

Then, we define the apparent homogenized matrix $A^*_{\eta}(x,\omega)$ by
\[
\forall 1 \leq i \leq d, \quad
A^*_{\eta}(x,\omega) e_i = \frac{1}{\eta^d} \int_{Y_{\eta}^x}
A_\epsilon(y,\omega) \nabla \chi_i(y,\omega) \, dy,
\]
where $e_i$ is the unit vector in the direction $i$ ($i=1,\dots,d$). We
denote this local homogenization procedure by $\mathcal{H}_\eta$, i.e.
\[
A_\eta^*(x,\omega)=\mathcal{H}_\eta (A_\epsilon(x,\omega)).
\]
This procedure is repeated at every macroscopic point 
(see Figure~\ref{rev} for illustration). Then, the coarse-scale
equation associated to~\eqref{microprob} is
\begin{equation}
\label{homprob}
-\text{div} (A^*_\eta(x,\omega) \nabla u^*)=f \text{ in $D$},
\end{equation}
with the same boundary conditions on $u^*$ as in~\eqref{microprob}. 

\subsubsection{Random microstructure and deterministic homogenized coefficients}

Homogenization of elliptic equations with random coefficients
has been extensively studied in the literature, and we refer
to~\cite{PV82,JKO94,blp,cd} for classical textbooks (see also the review
article~\cite{Bris-sg}). 
It is shown there that, if $\dis A_\epsilon(x,\omega) = 
A\left(\frac{x}{\epsilon},\omega\right)$ for some ergodic statistically
homogeneous (i.e. stationary) random field $A(x,\omega) \in \RR^{d
  \times d}$ (see 
e.g.~\cite{PV82,JKO94} for definitions), then the random solution
$u_\epsilon(\cdot,\omega)$ to~\eqref{microprob} 
converges, weakly in $H^1(D)$ and almost surely, to a deterministic
function $u^*$, solution to 
$$
-\text{div} (A^* \nabla u^*)=f \text{ in $D$},
$$
with appropriate boundary conditions (say $u^*=0$ on $\partial D$
if~\eqref{microprob} is complemented by $u_\epsilon(\cdot,\omega)=0$ on
$\partial D$). The homogenized coefficient, denoted $A^*$ in the
above equation, is a deterministic, constant matrix.

In addition, the numerical procedure outlined above is a practical way
to obtain a converging approximation of the homogenized matrix, in the
sense that
\begin{equation}
\label{hom112}
\lim_{\eta \rightarrow \infty} A^*_{\eta}(x,\omega) = A^*,
\end{equation}
almost surely, and for almost all $x$ (see~\cite{BP04}). Note that
(\ref{hom112}) can be equivalently written $\dis
\lim_{\epsilon\rightarrow 0} A^*_{\eta}(x,\omega)= A^*$ for any fixed
$\eta>0$. 

The only assumptions of ergodicity and stationarity do not allow for a precise 
convergence rate in~\eqref{hom112}. If, in addition, one assumes that
the matrix $A(x,\omega)$ decorrelates at large distances at some given
rate,
then one can also obtain a convergence rate in~\eqref{hom112} (see
e.g.~\cite{Y87,BP04}). A typical result is that
\begin{equation}
\label{eq:effcoef}
\EE \left[ \left| A^*_\eta(x,\cdot) - A^* \right|^2 \right] \leq 
C \left( \frac{\epsilon}{\eta}\right)^{\beta} \ \text{a.e.},
\end{equation}
for some $\beta>0$ and $C>0$ that depend on the decorrelation rate, but
are independent of $x$, $\eta$ and $\epsilon$, and where $\left| \cdot
\right|$ is any norm on the $d \times d$ matrices. 

Note that, in the absence of ergodicity, the homogenized coefficients
are a priori random matrices, that are invariant under the group of
actions representing homogeneous statistical fields.

\subsubsection{Stochastic homogenized coefficients}
\label{sec:stoch_homog}

As we mentioned in the introduction, the Multi-Level Monte Carlo method
is more efficient than a standard Monte Carlo method when the exact
homogenized coefficients are stochastic (otherwise, both methods are
equally efficient).
In stochastic homogenization, if no ergodicity is assumed, then the
homogenized coefficients can be stochastic.
In this work, we consider various cases in that setting.

The first case we consider is when the coefficient
in~\eqref{microprob} has the form
\[
A\left(x,{x\over\epsilon},\omega,\omega'\right) 
=
\widetilde{A}(x,\omega) \, B\left({x\over\epsilon},\omega'\right) \, \Id,
\]
where $\widetilde{A}$ and $B$ are two random scalar valued functions and
$\Id$ is the identity matrix. 
We thus see that $\omega$ corresponds to a randomness at the macroscopic
scale, while $\omega'$ corresponds to a randomness at the microscopic
scale. Let $A^*(x,\omega,\omega')$ be the homogenized matrix, which
depends on the macroscopic variables $(x,\omega)$, and also on the
microscopic randomness $\omega'$ as no ergodicity is assumed on $B$. 
We will assume that 
\[
\EE_{\omega'} \left[ \left| A^*(x,\omega,\omega')- \mathcal{H}_\eta \left(
      A\left(x,{x\over\epsilon},\omega,\omega'\right)\right) \right|^2 \right]
\leq C \left(\epsilon\over\eta \right)^\beta,
\]
where the constant
$C$ and the rate $\beta$ are independent of $\omega$, $x$, $\epsilon$
and $\eta$. 

A second, more general case we consider is when the randomness does
not explicitely split into a randomness at the macroscopic and the
microscopic scales. The heretogeneous field in~\eqref{microprob} then
writes $\dis A\left(x,{x\over\epsilon},\omega\right)$. We assume
that $A$ is scalar-valued, that we can do homogenization at
every macroscopic point, and that the following assumption holds: 
\[
\EE \left[ \left| A^*(x,\omega)-
  \mathcal{H}_\eta\left(A\left(x,{x\over\epsilon},\omega\right)\right)
\right|^2 \right] 
\leq C \left(\epsilon\over\eta \right)^\beta
\]
for some constant $C$ and rate $\beta$ independent of $x$, $\epsilon$
and $\eta$. This assumption is similar to the known results for ergodic
homogeneous stochastic homogenization recalled in~\eqref{eq:effcoef}. 

\subsection{Multi-Level Monte Carlo approach}
\label{sec:mlmc_gene}

We now briefly introduce the Multi-Level Monte Carlo (MLMC) approach in a
general context. The reader familiar with this approach can directly
proceed to Section~\ref{sec:meshes}. 

Let
$X(\omega)$ be a random variable. We are interested in the
efficient computation of the expectation of $X$, denoted by $\EE(X)$. In
our calculations below, $X$ is a function of the homogenized
coefficients or of the homogenized solutions. For example,
we are interested in the expectation of the homogenized coefficients
$\EE(A^*)$, or in the two-point covariance function. In this case, we
choose the random variable as 
$X(\omega) = \left[ A^*(x_1,\omega) \right]_{ij} \, 
\left[ A^*(x_2,\omega) \right]_{qp}$ for some $x_1 \in D$ and $x_2 \in
D$ (and some 
components $ij$ and $qp$ of the homogenized matrices). Other quantities
of interest include e.g. statistics of the homogenized solution. 

To compute an approximation of $\EE(X)$, a standard approach is the
Monte Carlo (MC) method. One first calculates a number $M$ of independent
realizations of the random variable $X$ (denoted $X^i$, $1 \leq i \leq
M$), and next approximates the expected value $\EE(X)$ by the arithmetic
mean (also called empirical estimator):
$$
E_M(X) := \frac{1}{M} \sum_{i=1}^{M} X^i.
$$

In this article, we are interested in Multi-Level Monte Carlo (MLMC)
methods. The idea is to consider the quantity of interest $X_l$ on
different levels $l$. In our case, levels denote various representative
volume sizes, or different mesh sizes. We assume that $L$ is the level
of interest, and that computing many realizations  at this level
is too computationally expensive.
We introduce levels smaller than $L$, namely $L-1,\dots,1$, and 
assume that the lower the level is, the cheaper the computation of $X_l$
is, and the less accurate $X_l$ is with respect to $X_L$. Setting $X_0=0$,
we write
$$
X_L = \sum_{l=1}^L \left( X_l - X_{l-1}\right).
$$
The standard MC approach consists in working with $M$ realizations of
the random variable $X_L$ at the level of interest $L$. In contrast,
within the MLMC approach, we work with $M_l$ realizations of
$X_l$ at each level $l$, with $M_1 \geq M_2 \geq \dots \geq M_L$.
We write
$$
\EE \left[ X_L \right] = \sum_{l=1}^L \EE \left[ X_l - X_{l-1} \right],
$$
and next approximate $\EE \left[ X_l - X_{l-1} \right]$ by an empirical
mean as above:
$$
\EE \left[ X_l - X_{l-1} \right] 
\approx
E_{M_l}(X_l - X_{l-1}) = \frac{1}{M_l} \sum_{i=1}^{M_l} \left( X^i_l -
  X^i_{l-1} \right),
$$
where $X^i_l$ is $i$th realization of the random variable $X$ computed
at the level $l$ (note that we have $M_l$ copies of $X_l$ and $X_{l-1}$,
since $M_l \leq M_{l-1}$). The MLMC
approach consists in approximating $\EE(X_L)$ by
\begin{equation}
\label{eq:approx_mlmc} 
E^L(X_L) := \sum_{l=1}^L E_{M_l}\left( X_l - X_{l-1}\right).
\end{equation}
As will be seen below, the realizations of $X_l$ used with those of
$X_{l-1}$ to evaluate $E_{M_l}\left( X_l - X_{l-1}\right)$
do not have to be independent of the realizations of $X_l$
used with those of $X_{l+1}$ to evaluate $E_{M_{l+1}}\left( X_{l+1} -
  X_l \right)$ (see also Remark~\ref{rem:iid} below).

In the following, we are interested in the root mean square errors
\begin{eqnarray}
\label{eq:e_mlmc}
e_{MLMC}(X_L) &=& 
\sqrt{\EE \left[ \left\| \EE(X_L)-E^L(X_L) \right\|^2 \right]},
\\
\label{eq:e_mc}
e_{MC}(X_L) &=& 
\sqrt{\EE \left[ \left\| \EE(X_L) - E_{M_L}(X_L) \right\|^2 \right]},
\end{eqnarray}
with an appropriate norm depending on the quantity of interest (e.g. the
absolute value for any entry of the homogenized coefficient, the $L^2(D)$
norm for the homogenized solution). 
For the error estimation, we will use (see e.g.~\cite{C98}) that, for
any random variable $X$, and any norm associated to a scalar product, 
\begin{equation}
\label{caf}
\EE \left[ \left\| \EE(X) - E_M(X) \right\|^2 \right] 
= 
\frac{1}{M} \EE \left[ \left\| X - \EE(X) \right\|^2 \right].
\end{equation}

\subsection{Definition of meshes and representative volume sizes}
\label{sec:meshes}

\begin{figure}[htp]
\center\scalebox{0.8}{\includegraphics{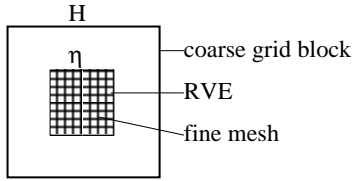}}
\caption{Parameters in the numerical homogenization procedure.}
\label{coarse}
\end{figure}

In our application, we will be dealing with various representative
volume sizes, 
and also possibly various sizes of coarse meshes
(see Figure~\ref{coarse} for illustration). In the framework of MLMC
approaches, choosing a level $l$ thus corresponds to choosing a
particular RVE size, \dots We denote the hierarchy of
coarse meshes on which we solve~\eqref{homprob} by 
\[
H_1 \geq H_2 \geq \dots \geq H_L.
\]
The number of realizations used at the level $i$ for
the coarse mesh size $H_i$ is denoted $M_i$. We take
\[
M_1 \geq M_2 \geq \dots \geq M_L.
\]
As for the representative volumes, we take their sizes according to
\[
\eta_1 \leq \eta_2 \leq \dots \leq \eta_L
\]
and the corresponding number of realizations is denoted 
\[
m_1 \geq m_2 \geq \dots \geq m_L.
\]
One could also use various fine-scale meshes for solving the local
representative volume problems~\eqref{eq:corrector}. We do not go in
this direction in this work. 

Note that the level $L$ always
corresponds to the most expensive choice (large RVE, or fine mesh), and
thus the smallest number of realizations. Note also that one does not
have to take the same number of levels $L$ for coarse-grid sizes and
RVEs.

\section{MLMC approach for the upscaled coefficients}
\label{sec:mlmc}

In this section, we describe how to use the MLMC approach to compute the
upscaled coefficients defined in Section~\ref{numHom} and the two-point
correlation functions. 
We focus on how to choose RVE sizes for the
problems~\eqref{eq:corrector}, and thus assume that these problems are
exactly solved. Setting 
$$
\delta_l(x)
=
\sqrt{\EE \left[ \left| A^*(x,\cdot) - A_l^*(x,\cdot) \right|^2 \right]},
$$
where $\left| \cdot \right|$ is some matrix norm,
we assume, following Section~\ref{sec:stoch_homog}, that
\begin{equation}
\label{eq:hyp}
\delta_l(x) \leq C \left( \frac{\epsilon}{\eta_l}
\right)^{\beta/2},
\end{equation}
for some $\beta > 0$ and $C>0$ independent of $l$, $\epsilon$, $\eta$
and of the macroscopic point $x \in D$ (in what follows, we keep the
dependency with respect to $x$ implicit in our notation).
For some special cases, one can obtain an estimate for $\beta$ rigorously.
For more complicated cases, we suggest in Section~\ref{sec:estim_beta}
below a pre-computation strategy that can
provide an estimate for $\beta$. 
Note that a Central Limit Theorem
type result corresponds to $\beta = d$ (see e.g.~\cite{BCLL11_mprf} for
such estimates in a weakly stochastic case). 

\medskip

For clarity, we summarize now our MLMC algorithm for
the upscaled coefficients:
\begin{enumerate}
\item Generate $m_1$ random variables $\omega_1$, \dots, $\omega_{m_1}$.
\item For each level $l$, $1\leq l \leq L$, and each realization
  $\omega_j$, $1\leq j \leq m_l \leq m_1$,  
\begin{itemize}
 \item Solve the RVE problems~\eqref{eq:corrector} on 
   $Y^x_{\eta_l}$: for any $i=1,\dots,d$,
$$
\text{div}(A_\epsilon(y,\omega_j)\nabla \chi_i(y,\omega_j)) = 0
\text{ in $Y^x_{\eta_l}$},
\quad
\chi_i(y,\omega_j) = y_i \ \text{on $\partial Y^x_{\eta_l}$}.
$$
\item Compute the homogenized matrix $A^*_l(x,\omega_j)$ with 
\[
\forall 1 \leq i \leq d, \quad
A^*_l(x,\omega_j) e_i
=
\frac{1}{{\eta_l}^d} \int_{Y^x_{\eta_l}} A_\epsilon(y,\omega_j) \, \nabla
\chi_i(y,\omega_j) \, dy.
\]
\end{itemize}
\item For each level $l$, $1\leq l \leq L$, compute
\[
E_{m_l}(A^*_l -A^*_{l-1}) = \frac{1}{m_l} \sum_{j=1}^{m_l} 
\left[ A^*_{l}(x,\omega_j) - A^*_{l-1}(x,\omega_j) \right], 
\ \
\text{with $A^*_0 =0$}.
\]
\item Compute the MLMC approximation $E^L(A_L^*)$ of the expected
  value $\EE(A_L^*(x,\cdot))$ following~\eqref{eq:approx_mlmc}:
\[
E^L(A_L^*) := \sum_{l=1}^L  E_{m_l}(A^*_l -A^*_{l-1}).
\]
\end{enumerate}

Let us now estimate the error in the approximation 
of $\EE \left( \left[ A^*_L \right]_{ij} \right)$, for any entry $ij$
($1 \leq i,j \leq d$) of the matrix $A^*_L$. To simplify the notation,
we write the calculations below as if $A^*_l$ were a scalar quantity
independent of $x$. These
calculations are to be understood as calculations on the entry $\left[
  A^*_l (x,\cdot) \right]_{ij} \in \RR$. 

For the MLMC approach, the error reads
\begin{equation*}
\begin{split}
e_{MLMC}(A^*_L) &= 
\sqrt{\EE \left[ \left( \EE(A^*_L) -  E^L(A_L^*) \right)^2 \right] } 
\\
&=  
\sqrt{\EE \left[ \left( \EE \left( \sum_{l=1}^L \left( A_l^* -A_{l-1}^*
        \right) \right) - \sum_{l=1}^L E_{m_l} \left( A_l^* -A_{l-1}^*
      \right) \right)^2 \right] }
\\
&= 
\sqrt{\EE \left[ \left( \sum_{l=1}^L \left( \EE -E_{m_l} \right) 
\left( A_l^*-A_{l-1}^* \right) \right)^2 \right] } 
\\
& \leq
\sum_{l=1}^L \sqrt{\EE \left[ \left( \left( \EE -E_{m_l} \right) 
\left( A_l^*-A_{l-1}^* \right) \right)^2 \right] } 
\\
&
=
\sum_{l=1}^L \frac{1}{\sqrt{m_l}} \sqrt{\EE \left[ \left( A_l^*-A_{l-1}^* - 
\EE(A_l^*-A_{l-1}^*) \right)^2 \right] } 
\end{split}
\end{equation*}
where we have used~\eqref{caf}.
Writing that $A_l^*-A_{l-1}^* = (A_l^*-A^*) + (A^*-A_{l-1}^*)$, and
since $m_l \leq m_{l-1}$, we deduce that
\begin{equation*}
\begin{split}
e_{MLMC}(A^*_L) 
&
\leq 
\frac{1}{\sqrt{m_L}}\sqrt{ \EE \left[ \left( A_L^*-A^*- \EE
      \left(A_L^*-A^* \right)\right)^2 \right]} 
\\
& +
\sum_{l=1}^{L-1} \frac{2}{\sqrt{m_{l+1}}} \sqrt{ \EE \left[ \left( A_l^*-A^*-
    \EE \left(A_l^*-A^* \right) \right)^2 \right] } + 
\frac{1}{\sqrt{m_1}} \sqrt{\EE \left[ (A^*)^2 \right] }
\\
&
\leq \sum_{l=1}^L \frac{2}{\sqrt{m_{l+1}}} 
\sqrt{\EE \left[ \left( A_l^*-A^* \right)^2 \right] }
+ \frac{1}{\sqrt{m_1}} \sqrt{\EE \left[ (A^*)^2 \right] }
\\
&
\leq \sum_{l=1}^L \frac{2}{\sqrt{m_{l+1}}} \delta_l 
+ \frac{1}{\sqrt{m_1}} \sqrt{ \EE \left[ (A^*)^2 \right]},
\end{split}
\end{equation*}
where, for ease of notation, we have introduced some $m_{L+1} \leq 
m_L$. Using~\eqref{eq:hyp}, we deduce that
$$
e_{MLMC}(A^*_L) 
\leq 
C \sum_{l=1}^L \frac{1}{\sqrt{m_{l+1}}} 
\left( \frac{\epsilon}{\eta_l} \right)^{\beta/2}
+ \frac{1}{\sqrt{m_1}} \sqrt{ \EE \left[ (A^*)^2 \right]}.
$$

For a fixed error, 
the optimal choice for the number $m_l$ of realizations at level $l$
(namely for the RVE of size $\eta_l$) is reached when these error parts are
equilibrated. Therefore, we choose 
\begin{equation}
\label{eq:m_l}
m_l = 
\begin{cases}
\dis
\left( \frac{\eta_L}{\epsilon}\right)^\beta \EE \left[ (A^*)^2 \right] 
\alpha_1^{-2}, &l=1,
\\
\dis
\left( \frac{\eta_L}{\eta_{l-1}}\right)^\beta \alpha_l^{-2}, & 2 \leq l \leq L+1,
\end{cases}
\end{equation}
for some parameters $\alpha_l$, and we check that indeed $m_{L+1} \leq
m_L$, provided $\alpha_{L+1} = \alpha_L$. 
We then have
\begin{equation}
e_{MLMC}(A^*_L) 
\leq
C \left( \frac{\epsilon}{\eta_L} \right)^{\beta/2} \sum_{l=1}^L
\alpha_l.
\label{eq:stat_error_mlmc}
\end{equation}

\medskip

For comparison, we consider the error if we calculate the approximated
upscaled coefficient only for the largest RVE (of size $\eta_L$), using
a standard MC method with $\widehat{m}_L$ independent
samples. Using~\eqref{caf}, we find that the MC error reads
$$
e_{MC}(A^*_L)
= 
\sqrt{ \EE \left[ \left( \EE \left(A^*_L\right) - 
E_{\widehat{m}_L}\left(A^*_L\right) \right)^2 \right] }
=
\frac{1}{\sqrt{\widehat{m}_L}} \sqrt{ \EE \left[ 
\left( A^*_L - \EE\left(A^*_L\right) \right)^2 \right]}.
$$
As pointed out above, $A^*$ is assumed to be a random quantity, with
some positive variance. It is thus natural to assume that the variance
of $A^*_L$ is roughly independent of $L$, and hence that the MC
error is of the order of $C / \sqrt{\widehat{m}_L}$. 
To have an error of the same order as that given by the MLMC approach,
we take $\dis \widehat{m}_L = O \left( \left(
    \frac{\eta_L}{\epsilon} \right)^\beta \right)$ independent
realizations. 

\medskip

\begin{figure}[htp]
\center
\subfigure[$\dis \epsilon=\frac{\eta_1}{10}$]{
\scalebox{.4}{ \includegraphics{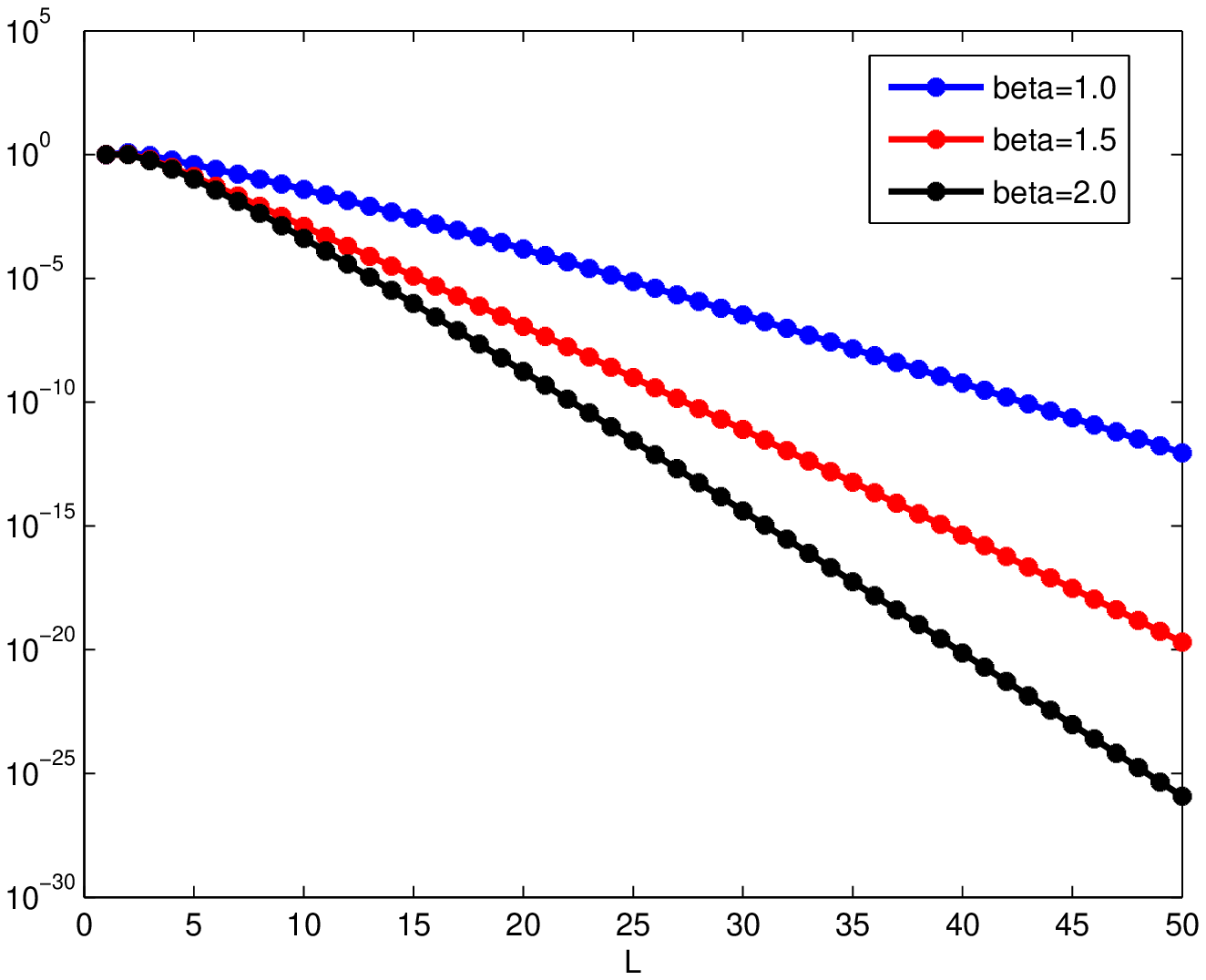}}}
\subfigure[$\dis \epsilon=\frac{2^{-50}}{10}$]{
\scalebox{.4}{ \includegraphics{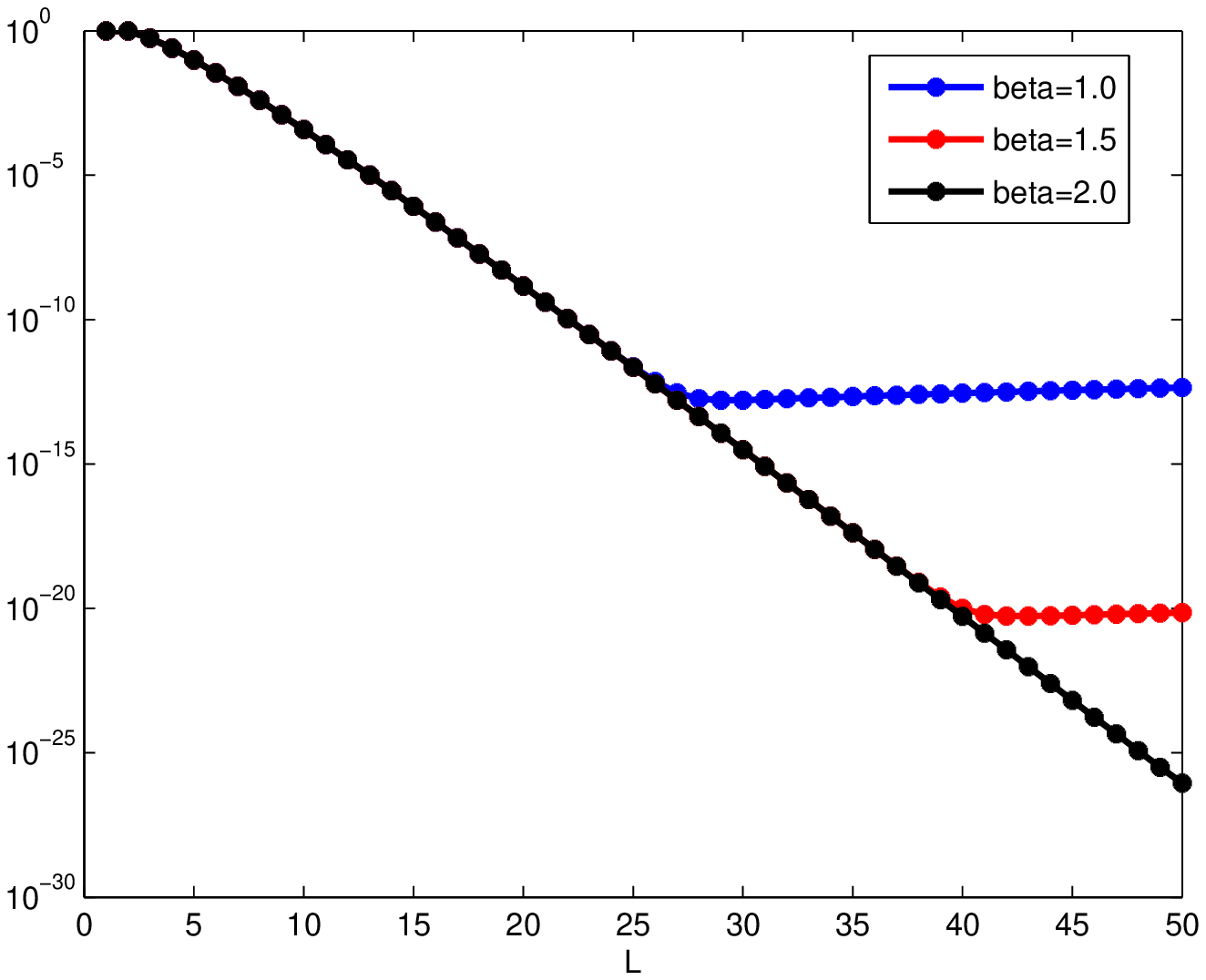}}}
\caption{RVE cost ratio 
$\dis \frac{W^{MLMC}_{\text{RVE}}}{W^{MC}_{\text{RVE}}}$ for different
numbers of levels $L$ and rates $\beta$ (we work with  
$\eta_l=2^{l-L}$, $\alpha_l=1/L$ for all $l$, $d=2$ and
$\EE\left[(A^*)^2\right]=1$).}
\label{rvework}
\end{figure}

Now that we have chosen the number of realizations for both approaches so
that they reach the same accuracy, we
are in position to compare their cost. Let $N_l$ denote the cost of
solving the RVE problem~\eqref{eq:corrector} on the domain
$Y^x_{\eta_l}$ of size $\eta_l^d$. 
The number of degrees of freedom needed is of the order of
$(\eta_l/\epsilon)^d$. Assuming that $N_l = (\eta_l/\epsilon)^d$, the
MLMC cost is $\dis W^{MLMC}_{\text{RVE}} = \sum_{l=1}^L m_l N_l$, hence 
$$
W^{MLMC}_{\text{RVE}}
=
\sum_{l=2}^L \left( \frac{\eta_L}{\eta_{l-1}}\right)^\beta \alpha_l^{-2} 
\left(\frac{\eta_l}{\epsilon}\right)^d 
+
\left( \frac{\eta_L}{\epsilon}\right)^\beta \EE \left[ (A^*)^2 \right]
\alpha_1^{-2} \left(\frac{\eta_1}{\epsilon}\right)^d.
$$
In the case of the MC approach, the cost reads
$$
W^{MC}_{\text{RVE}} 
= 
\widehat{m}_L N_L 
=
C \left( \frac{\eta_L}{\epsilon} \right)^\beta \left(
  \frac{\eta_L}{\epsilon} \right)^d
=
C \left( \frac{\eta_L}{\epsilon}\right)^{\beta+d}.
$$
On Figure~\ref{rvework}, we plot the ratio 
$\dis \frac{W^{MLMC}_{\text{RVE}}}{W^{MC}_{\text{RVE}}}$ for different
numbers of levels $L$ and rates $\beta$, with the choice
$\eta_l=2^{l-L}$. Note then that the largest RVE is always of size
$\eta_L = 1$, independently of $L$, and that the smallest RVE size depends on
$L$, and is $\eta_1=2^{1-L}$. On the right plot, we consider the case
when $\epsilon$ is fixed at a very small value
independent of $L$. This value is sufficiently
small to ensure that, even for the largest considered $L$, the smallest
RVE is larger than $\epsilon$
(thereby ensuring scale separation). On the left plot, we consider
a more practical situation (which is the regime we choose for our
numerical experiments of 
Section~\ref{sec:num}), when $\epsilon$ depends on $L$ and is always 10
times smaller that the smallest RVE. This
leads to values of $\epsilon$ that are larger (and thus easier to handle
numerically) than that considered on the right plot.
 
As we can see, for a given number $L$ of levels, the larger the rate
$\beta$ is, the smaller the cost ratio is, at equal accuracy. Otherwise
stated, the faster the convergence of the apparent homogenized matrix 
with respect to the RVE size, the more efficient the MLMC approach
is. We also observe on the right plot that, at fixed $\beta$ and
$\epsilon$, the gain in terms of cost first increases when $L$ increases
and then reaches a plateau for large $L$.  

\begin{bem}
In the above calculations, we have assumed that the cost of solving
a local problem scales {\em linearly} with the number ${\cal N}$ of degrees
of freedom. This is true if one uses iterative solvers 
and the condition number of the preconditioned system
is independent of the small scale $\epsilon$. One can also compare the cost
between the MLMC and MC approaches under different assumptions
(e.g. when the cost of solving a local problem scales as $C(\epsilon)
{\cal N}^{1+\gamma}$ for some $\gamma \geq 0$).
\end{bem}

\begin{bem}
To compute the optimal number of realizations following~\eqref{eq:m_l},
one needs to know the value $\beta$ of the rate
in~\eqref{eq:hyp}. In general, this rate is not
analytically known.
To address this difficulty, we propose in
Section~\ref{sec:estim_beta} below some means to estimate the value of
$\beta$ based on a priori, offline computations. 
\end{bem}

We have shown above how to estimate $\EE(A^*_L(x,\cdot))$ at any macroscopic
point $x$. Another important quantity is the two-point correlation
function 
$$
\Cor_{A^*}(x,y) := \EE \left( \left[ A^*(x,\omega) \right]_{ij} \
\left[ A^*(y,\omega) \right]_{qp} \right)
$$
between the components $ij$ and $qp$ of the homogenized matrix at points
$x$ and $y$ (note that we work with non-centered values of $A^*$).
For simplicity, we only consider two fixed locations $x$ and
$y$. Consider $m_l$ independent realizations of the homogenized matrix
$A^{*,k}_l$ at level $l$ ($1 \leq k \leq m_l$). We define
$$
Cor_{m_l}(A^*_l) := \frac{1}{m_l} \sum_{k=1}^{m_l} \left[ A^{*,k}_l(x)
\right]_{ij} \ \left[ A^{*,k}_l(y) \right]_{qp}
$$
as an empirical estimator of $\EE \left( \left[ A^*_l(x,\omega) \right]_{ij} \
\left[ A^*_l(y,\omega) \right]_{qp} \right)$.
The MLMC approximation of the two-point correlation function 
$\Cor_{A^*}(x,y)$ then reads
$$
Cor^L(A^*_L) := \sum_{l=1}^L \left(Cor_{m_l}(A^*_l)-Cor_{
    m_l}(A^*_{l-1})\right).
$$

\begin{bem}
We have considered above that we could exactly solve
the RVE problems~\eqref{eq:corrector}. In practice, these problems are
solved numerically, within some accuracy. 
A natural extension of assumption~\eqref{eq:hyp} is to assume that
the error in the approximation of $A^*$ (due to working on a truncated
domain $Y_{\eta_i}$ of size $\eta_i$ with a finite discretization on a
mesh of size $h_j$) satisfies 
\[
\delta_{ij} \leq C \sqrt{\left({\epsilon\over \eta_i}  \right)^\beta + 
 \left({h_j\over \epsilon}  \right)^\gamma},
\]
for some constant $C$ independent of $h_j$, $\eta_i$ and $\epsilon$. 
An analysis similar to the one above then follows. Note also that it may
be possible to solve the local problems on some RVEs with a coarser
approximation and correct this using the nearby RVEs, computed at full
accuracy, in the spirit of the strategy proposed in~\cite{brown} in
another context. The adaptation of such an idea to our context goes
beyond the scope of the current work. 
\end{bem}

\section{MLMC for the homogenized solution}
\label{sec:homog_sol}

In this section, we show how to estimate the expectation of the
homogenized solution using the MLMC approach. We also introduce an
extension of that approach, namely the weighted MLMC approach, in
Section~\ref{solutionGeneral}.

\subsection{Separable case}
\label{sec:separable}

In this section, we assume that the coefficient in~\eqref{microprob} reads 
$$
A_\epsilon(x,\omega,\omega') = \widetilde{A}(x,\omega) \, 
B\left(\frac{x}{\epsilon},\omega'\right) \, \Id
$$
for two scalar valued functions $\widetilde{A}$ and $B$,
and therefore satisfies a separation of scales assumption. 
The coarse-scale problem associated to the highly oscillatory
problem~\eqref{microprob} is
$$
-\text{div} \left[ \widetilde{A}(x,\omega) \, B^*(\omega') \nabla u^* \right]
= f \text{ in $D$}.
$$
We expect most of the randomness of the coefficient at the
coarse-scale to be in $\widetilde{A}(x,\omega)$. We thus use a simplistic
treatment for averaging over $\omega'$ and approximate the above
equation by
\begin{equation}
\label{eq:coarse-PDE}
-\text{div} \left[ \widetilde{A}(x,\omega) \, 
\EE_{\omega'} \left[ B^* \right] \nabla u^* \right]
= f \text{ in $D$}.
\end{equation}
We are going to compute an approximation of $\EE(u^*)$, using 
the tuples $\left(H_l, M_l, \eta_l,
  m_l\right)$ for $1\leq l\leq L$. 

We first need to calculate the
homogenized coefficient $B^*(\omega')$. To do so, we solve in each
direction, $1 \leq j\leq d$, and for each realization 
$\dis B^i\left(\frac{y}{\epsilon},\omega' \right)$ of the coefficient, 
$1 \leq i \leq m_l$, the RVE problem
\begin{equation}
\label{eq:RVE_sep}
\begin{array}{rcl}
\dis -\text{div} \left[ B^i\left(\frac{y}{\epsilon},\omega' \right)
\nabla \chi^i_j \right]
&=& 0 \ \ \text{in $Y_{\eta_l}$},
\\ \noalign{\vskip 3pt}
\chi^i_j(y,\omega') &=& e_j \cdot y  \text{ on $\partial Y_{\eta_l}$},
\end{array}
\end{equation}
and calculate the corresponding homogenized coefficient:
$$
\left[ B^{*i}_l\right]_{n,m}
= \frac{1}{|Y_{\eta_l}|}
\int_{Y_{\eta_l}} \nabla \chi^i_n \cdot B^i \nabla \chi^i_m.
$$
Note that we have kept implicit the dependency of $\chi^i_j$ with
respect to the level $l$. We then introduce
$$
E_{m_l}(B_l^*)= \frac{1}{m_l} \sum_{i=1}^{m_l} B^{*i}_l,
$$
which is an approximation (at level $l$) of 
$\EE_{\omega'} \left[ B^* \right]$. We correspondingly introduce $u_l$,
solution to 
$$
-\text{div} \left[ \widetilde{A}(x,\omega) \, E_{m_l}(B_l^*) \nabla u_l
\right] = f \text{ in $D$}.
$$
In turn, this equation is solved on a mesh of size $H_l$, for several
realizations of $\widetilde{A}(x,\omega)$.
We thus eventually define $u^k_l$ (with $1 \leq l \leq L$
and $1 \leq k \leq M_l$),
solution (on a mesh of size $H_l$) to the coarse-scale equation 
\begin{equation}
\label{eq:homsoln}
-\text{div} \left[ \widetilde{A}^k(x,\omega) \, E_{m_l}(B_l^*) \nabla u^k_l
\right] = f \ \text{ in $D$}.
\end{equation}
The expected value $\EE(u_l)$ is approximated in a
standard Monte Carlo fashion by 
$$
\EE(u_l) 
\approx 
E_{M_l}(u_l) := \frac{1}{M_l} \sum_{k=1}^{M_l} u^k_l,
$$
where $u^k_l$ is the solution to (\ref{eq:homsoln}). 

\medskip

To approximate our quantity of interest, $\EE(u_L)$, we can first
perform the above procedure only at the level $L$. This yields a
standard Monte Carlo approximation of $\EE(u_L)$. 

An alternative approximation is that provided by the MLMC approach,
which reads
$$
E_L(u_L) := \sum_{l=1}^L E_{M_l}(u_l-u_{l-1})
\quad \text{with} \quad
u_0 = 0.
$$
Introducing the norm
$$
\left\| X \right\| = 
\left[ \EE \left( \left\| X(\omega) \right\|_{H^1(D)}^2 \right) \right]^{1/2},
$$
the MLMC error is estimated following the same lines as in
Section~\ref{sec:mlmc}. We obtain
$$
\left\| \EE(u_L) - E_L(u_L) \right\|
\lesssim
\sum_{l=1}^L \frac{1}{\sqrt{M_l}} \left\| u^*-u_l \right\|
+ \frac{1}{\sqrt{M_1}} \left\| u^* \right\|.
$$
To bound from above $\left\| u^*-u_l \right\|$, we introduce
$u^*_{H_l}$, approximate solution to~\eqref{eq:coarse-PDE} on a mesh of size
$H_l$. It follows that 
$$
\left\| u^*-u_l \right\|
\leq
\left\| u^* - u^*_{H_l} \right\| + \left\| u^*_{H_l} - u_l \right\|.
$$
The first term is a discretization error, which typically satisfies
(e.g. if we use a P1 Finite Element method) the bound
$\left\| u^* - u^*_{H_l} \right\| \lesssim H_l$. 
For the second term, it holds (all expectations are taken w.r.t. $\omega'$)
\begin{eqnarray*}
\left\| u^*_{H_l} - u_l \right\|
& \lesssim & 
\sqrt{\EE\left[ \left| E_{m_l}(B_l^*) - \EE \left( B^* \right) \right|^2 \right]}
\\
& \leq &
\sqrt{\EE \left[ \left| E_{m_l}(B_l^*)- \EE(B^*_l) \right|^2 \right]} 
+  
\left| \EE(B^*_l) - \EE \left( B^* \right) \right|
\\
&\lesssim& 
\frac{1}{\sqrt{m_l}} \sqrt{\EE\left[ \left| B_l^*- \EE(B^*_l) \right|^2 \right]}
 + \delta_l.
\end{eqnarray*}
Using our assumption~\eqref{eq:hyp}, that is 
$\dis
\delta_l^2 \lesssim \left( \frac{\epsilon}{\eta_l} \right)^\beta$
for some $\beta > 0$, and assuming that the variance of $B^*_l$ is
essentially independent of $l$, we get
$$
\left\| \EE(u_L) - E_L(u_L) \right\|
\lesssim 
\sum_{l=1}^L \frac{1}{\sqrt{M_l}} \left( H_l +
\left( \frac{\epsilon}{\eta_l} \right)^{\beta/2}
+ \frac{C}{\sqrt{m_l}} \right) +\frac{1}{\sqrt{M_1}}.
$$
For the standard MC approach (with $\widehat{M}$ independent samples),
the error reads 
$$
\left\| \EE(u_{MC})- E_{\widehat{M}}(u_{MC}) \right\|
\lesssim
\frac{1}{\sqrt{\widehat{M}}},
$$
provided the variance of $u_{MC}$ is of order one.

\subsection{General, non-separable case}
\label{solutionGeneral}

In general, the coefficient in~\eqref{microprob} is of the form $\dis
A\left(x, \omega, \frac{x}{\epsilon}\right)$, where there is no
separation between the 
macroscopic and the microscopic randomness. In this case, the RVE
problems are parameterized by the macroscale position $x$, and thus need
to be solved in each coarse-grid block (in contrast to the separable
case considered in Section~\ref{sec:separable}, where the local RVE
problem~\eqref{eq:RVE_sep} is independent of $x$).

At any level $l$, let $N_l \propto H_l^{-d}$ be the number of
coarse-grid blocks. We denote by $\mathcal P_l$ the set of the
macroscale grid points at which we solve a RVE problem, with
$\text{Card } \mathcal P_l = N_l$. We assume that the coarse grids are
nested from one level to the other, so that 
$\mathcal P_1 \subset \mathcal P_2 \subset \cdots \subset \mathcal P_L$. 
As before, on each grid of size $H_l$, we solve $M_l$ coarse grid
problems. 
To calculate the effective coefficient, we solve the RVE problems at
each coarse grid point and for each realization of $A$, and we next
average the energy over the spatial domain. 
Since the sets $\mathcal P_l$ are nested, once we have computed
$A^*_{\eta_l,H_l}$ (using RVEs of size $\eta_l$) at the macroscopic
points of the coarse mesh of
size $H_l$, we readily get $A^*_{\eta_l,H_j}$ for $j<l$ (see
Table~\ref{illus}). Thus, at each level $l < L$, and at each point of the
grid of mesh size $H_l$, we only have to solve 
$M_l-M_{l+1}$ RVE problems (associated to independent realizations) on
RVEs of size $\eta_l$, and not $M_l$ of them. 

\begin{table}[htbp]
\center
\begin{tabular}{c|cccc|c}
             & $H_L$ & $H_{L-1}$ & $\cdots$ & $H_1$ & \# coefficients to calculate\\
             &       &          &         &       &  with RVE size $\eta_l$\\
\hline
$\eta_1$     &       &           &         &  \textcolor{blue}{$A^*_{\eta_{1}, H_{1}}$ }     & $M_1-M_2$\\
$\vdots$    &       &           &        &  $\vdots$     & $\vdots$\\
$\eta_{L-1}$ &       & \textcolor{blue}{ $A^*_{\eta_{L-1}, H_{L-1}}$ } & $\cdots$ & $A^*_{\eta_{L-1}, H_{1}}$   & $M_{L-1}-M_L$\\
$\eta_{L}$   &  \textcolor{blue}{ $A^*_{\eta_L, H_L}$ }    &  $A^*_{\eta_L, H_{L-1}}$         &   $\cdots$      &   $A^*_{\eta_L, H_{1}}$      & $M_L$\\
\noalign{\vskip 5pt} 
\hline
\# coefficients & & & & & \\
 on grid size $H_l$  & $M_L$ & $M_{L-1}$ &$\cdots$ & $M_1$ &  \\
\end{tabular}
\caption{Calculating the coefficients on the diagonal (shown in blue)
  will automatically give the lower triangular values in the
  matrix.}
\label{illus} 
\end{table}

We denote $u^*_{\eta_j, H_i}$ the solution to the coarse-scale equation
discretized on a grid of size $H_i$, and where the effective coefficient
is computed from local problems set on RVEs of size $\eta_j$.

To approximate $\EE(u^*_{\eta_L, H_L})$, we can use a MLMC approach
based on the
solutions $u^*_{\eta_j, H_j}$, $1 \leq j \leq L$. However, such an
approach discards the solutions $u^*_{\eta_j, H_i}$ for $i < j$,
which are however easy to compute. Indeed, once the coefficient
$A^*_{\eta_j, H_j}$ has been obtained at some level $j$, computing the
solutions $u^*_{\eta_j, H_i}$ for {\em all} meshes $i \leq j$ is as
inexpensive as computing $u^*_{\eta_j, H_j}$ only for the mesh $j$.

To benefit from this fact, we can approximate $\EE(u^*_{\eta_L, H_L})$ using
a {\em weighted} MLMC approach, which is defined as
\begin{equation}
\label{eq:w-mlmc_approx}
E^{L*}_{weighted} := \sum_{l=1}^L \alpha_l E^*_{M_l}(u^*_l-u^*_{l-1}),
\end{equation}
where $\alpha_l$ ($1 \leq l \leq L$) are parameters to be determined, and
\begin{eqnarray*}
E^*_{M_l}(u^*_l-u^*_{l-1})
& := &
\frac{1}{M_l} \sum_{j=l}^L \left( \sum_{i=1}^{M_j-M_{j+1}}
\left( u^*_{\eta_j, H_l} - u^*_{\eta_j, H_{l-1}} \right)(\omega_i) \right)
\\
&=& 
\frac{1}{M_l} \sum_{j=l}^L (M_j-M_{j+1}) 
E_{M_j-M_{j+1}}(u^*_{\eta_j, H_l} - u^*_{\eta_j, H_{l-1}}),
\end{eqnarray*}
where we have set $M_{L+1}=0$ and $u_{\eta_j,H_0}=0$ for any $1\leq j
\leq L$. Note that if $\alpha_l = 1$ for all $l$, we recover the
standard MLMC approach.

Errors associated to the weighted MLMC approach are estimated in
Appendix~\ref{sec:app}. 


\section{Numerical results}
\label{sec:num}

We consider the problem 
$$
-\text{div}\left[ A\left(x, \omega, \frac{x}{\epsilon}, \omega'\right)
  \nabla u_\epsilon \right] = f \text{ in $D = (0,1)^d$},
$$
complemented by boundary conditions that will be made precise
below. Likewise, the function $f$ will be given below. Note that 
the exact homogenized coefficient is independent of these choices. 

In what follows, we compare our MLMC results with standard MC results at
the highest level. We equate the cost for calculating
the coefficient and the solution separately and compare the errors (in
contrast to the theoretical analysis of Sections~\ref{sec:mlmc}
and~\ref{sec:homog_sol}, where we have equated the accuracies and
compared the costs). 

We will consider both one-dimensional and two-dimensional examples.
For the one-dimensional cases, we have implemented the method in Matlab,
and used the analytical solutions of the various PDEs. In the
two-dimensional cases, we use a rectangular mesh with cell-centered 
finite volumes. To solve the PDEs, we use the modular toolbox DUNE, the
Distributed and Unified Numerics Environment~\cite{dunepaperI:08,
  dunepaperII:08, dune:Fem, BlattBastian2007}. 

When we use the MLMC approach to approximate the homogenized
coefficient, we consider $L=3$ different RVE sizes $\eta_l$, 
unless specified otherwise. Likewise,
when we compute the 
homogenized solution, we also use $L=3$ different coarse grids of mesh
size $H_l$. For all the computations, we have used the same fine grid
(see Table~\ref{data}). We have made sure that the smallest RVE we
consider is much larger than the characteristic length scale $\epsilon$
(given below for each example)
of the field $A$. 

\begin{table}[htbp]
\center
 \begin{tabular}{ccclc}
$l$ & $H_l$           & $h_l$            &$\eta_l$ & \# cells in RVE of size $\eta_l$\\
\toprule\\
$1$ & $\dfrac{1}{16}$ & $\dfrac{1}{128}$ & $0.125 $ & $256$\\[1.5ex]
$2$ & $\dfrac{1}{32}$ & $\dfrac{1}{128}$ & $0.25 $  & $1024$\\[1.5ex]
$3$ & $\dfrac{1}{64}$ & $\dfrac{1}{128}$ & $0.5 $   & $4096$\\
 \end{tabular}
\caption{Parameters for the MLMC approach (two-dimensional cases).}
\label{data}
\end{table}

We explain in Section~\ref{sec:estim_beta} how to numerically estimate 
the rate of convergence $\beta$ in~\eqref{eq:hyp}. In
Section~\ref{sec:num_coeff}, 
we present numerical results for the homogenized coefficients.
Next, in Section~\ref{sec:num_sol}, we present numerical results for
homogenized solutions.

\subsection{Numerical study of the convergence rate}
\label{sec:estim_beta}

In our theoretical study described above, we have assumed that 
$$
\EE \left[ \left| A^*_l - A^* \right|^2 \right] 
\leq 
C \left( \frac{\epsilon}{\eta_{l}}\right)^\beta
$$
for some constant $C$ and rate $\beta$ independent of $\epsilon$ and
$\eta$ (see~\eqref{eq:hyp}). In this section, we numerically estimate
the parameter $\beta$ on a practical example. 

The considered scalar coefficient $\dis A
\left(\frac{x}{\epsilon}, \omega' \right)$ (defined for $x \in D \subset
\RR^2$) is a random field with expected value $\EE(A)=10$ (independent
of $x$ and $\epsilon$) and Gaussian covariance function 
$$
\text{cov}(x,x') = \Cov \left[ A \left( \frac{x}{\epsilon}, \cdot \right),
A \left( \frac{x'}{\epsilon}, \cdot \right) \right]
=
\sigma^2 \exp \left( -\frac{ | x-x' |^2}{\epsilon^2 \tau_0^2} \right),
$$ 
with $\sigma=\sqrt{2}$, $\tau_0 = \sqrt{2}$ and $\tau = \epsilon \tau_0 =
0.04$ (recall that $| x-x' |$ denotes the Euclidean distance in $\RR^2$). We
generate samples of the coefficient with the Karhunen-Lo\`eve
expansion. By construction, the characteristic length scale $\epsilon$
is related to the correlation length in $\text{cov}(x,x')$, which is of
the order of $\epsilon \tau_0$.

\begin{figure}\center
\psfrag{ln\(epsilon/eta}{\LARGE $\ln \epsilon/\eta_l$}
\psfrag{l}{}
\psfrag{\)}{}
\scalebox{.4}{\includegraphics{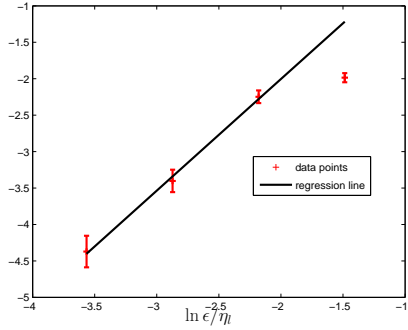}}
\caption{Using~\eqref{eq:estim_beta} to estimate $\beta$: computed data
  points (along with error bars) and the 
  corresponding linear regression line with slope $\beta=1.53$.}
\label{betareg} 
\end{figure}

For any $1\leq l \leq L$, we calculate the effective coefficients
$A^*_l(\omega'_j)$ for the RVE $[0, \eta_l]^2$ (with $\eta_l=
0.5^{L-l}$) for various realizations $\omega'_j$, $1 \leq j \leq
m_l$. The theoretical reference value is 
$\dis A^* = \lim_{\eta \to \infty} \EE \left( A^*_\eta \right)$, to which
we cannot access in practice. We thus define the reference value
as
$$
A^*_{ref} := \frac{1}{L} 
\sum_{l=1}^{L} \frac{1}{m_l} \sum_{j=1}^{m_l} A^*_l(\omega'_j),
$$
where we have taken into account all the realizations on the RVEs $[0,
\eta_l]^2$, $1\leq l \leq L$, in order to decrease as much as possible
the statistical error.
In practice, we work with $L=4$ and 
$\mathfrak{m} = (m_1,m_2,m_3,m_4)= (2000, 1000, 300, 140)$.
For each RVE of size $\eta_l$, $1\leq l\leq L$, we expect
from~\eqref{eq:hyp} that 
$$
\frac{1}{m_l}\sum_{j=1}^{m_l} \left| A^*_l(\omega'_j)- A^*_{ref} \right|^2
\approx
\EE \left[ \left| 
A^*_l - A^*_{ref}
\right|^2 \right]
\approx
C \left( \frac{\epsilon}{\eta_l}\right)^{\beta},
$$
hence
\begin{equation}
\label{eq:estim_beta}
\ln \left( 
\frac{1}{m_l}\sum_{j=1}^{m_l} \left| A^*_l(\omega'_j)- A^*_{ref} \right|^2
\right) \approx \beta \ln \left( \frac{\epsilon}{\eta_l}\right) + \ln C.
\end{equation}
Results are shown on Figure~\ref{betareg}, where we plot the computed
data points (with error bars) and the corresponding linear regression line. 
We see that we find a straight line with slope $\beta=1.53$ and intercept
$\ln C=1.059$ in the asymptotic regime $\eta \gg 1$. Note that the value
of $\beta$ is smaller than, but close 
to, the value $\beta^{\rm theo} = d = 2$ that would be obtained using a Central
Limit theorem argument (see discussion below~\eqref{eq:hyp}). 
In the numerical tests that follow, we will often consider only the
three smallest RVE ($\eta = 0.5$, 0.25 and 0.125), for computational
cost reasons. The slope of the regression line computed on the
basis of these three smallest RVE decreases to $\beta = 1.0095$. 

These estimations will be useful in Section~\ref{sec:coeff_2D} below
(see Example 1).

\subsection{Computation of the homogenized coefficient}
\label{sec:num_coeff}

We first consider the one-dimensional situation
(Section~\ref{sec:coeff_1D}) and next turn to two-dimensional test cases
in Section~\ref{sec:coeff_2D}. 

\subsubsection{One dimensional examples}
\label{sec:coeff_1D}

Since the local problems~\eqref{eq:corrector} are analytically solvable, 
we can afford to take many levels and many realizations at each
level.

\paragraph{Example 1 (separable coefficient)}

As a first test-case, we consider a coefficient $\dis
A\left(\frac{x}{\epsilon},\omega, \omega'\right)$ such that its inverse
reads 
$$
A^{-1}\left(\frac{x}{\epsilon}, \omega, \omega'\right)
=
\left[ C  +\sum_{i=1}^N \chi_i(\omega') \sin^2\left(\frac{2 \pi
      x\varphi_i}{\epsilon}\right)\right] \exp(\omega),
$$
where $\omega$ and $\chi_i$ are i.i.d. random variables, uniformly
distributed in $[0,1]$, $\varphi_i$ are fixed random numbers in
$[0.2,2]$, and $C>0$ is a deterministic constant. Note that $A^{-1}$ is
uniformly bounded away from 0. This coefficient is {\em separable} in the
sense that $A^{-1}$ writes as a {\em product} of
a function of $\omega$ times a function of $\omega'$. 
For a fixed realization $\omega$, it is well known that, in the
one-dimensional situation, the homogenized
coefficient is the harmonic mean. Therefore the apparent homogenized
coefficient on the RVE $[a,b]$ is
\begin{eqnarray*}
& & A^*_{a,b}(\omega, \omega')\\
&=&\left( \frac{1}{b-a} \int_a^b A^{-1}\left(\frac{x}{\epsilon}, \omega, \omega'\right) \right)^{-1}\\
&=&\left( \frac{\exp(\omega)}{b-a}\left[ C(b-a)
+\sum_{i=1}^N \chi_i (\omega')\left( \frac{b-a}{2}- 
\frac{\sin\left(4 \pi b \varphi_i/\epsilon\right) - \sin\left(4 \pi a \varphi_i/\epsilon\right)}{ (8 \pi \varphi_i)/\epsilon} 
\right)\right]\right)^{-1}.
\end{eqnarray*}
In our simulation, we use the values $C=1$, $N=20$ and $\dis
\epsilon=\frac{0.5^L}{10}$ (which ensures that the smallest RVE
considered in the MLMC approach, of size $0.5^L$,
is much larger than $\epsilon$, the characteristic length of
the heterogeneous coefficient).
As reference, we use the MC approach with $1000$ realizations of the
apparent coefficient on the largest RVE $[a_L,b_L] = [0,0.5]$. 
In what follows, a realization is determined by the tuple $(\omega,
\chi_1(\omega'),\dots,\chi_N(\omega'))$. Likewise, expectations are
taken with respect to $\omega$ and $\omega'$.

\begin{figure}[htp]\center
\subfigure[Relative mean square errors on the expected value of the effective coefficient.]
{\scalebox{.4}{ \includegraphics{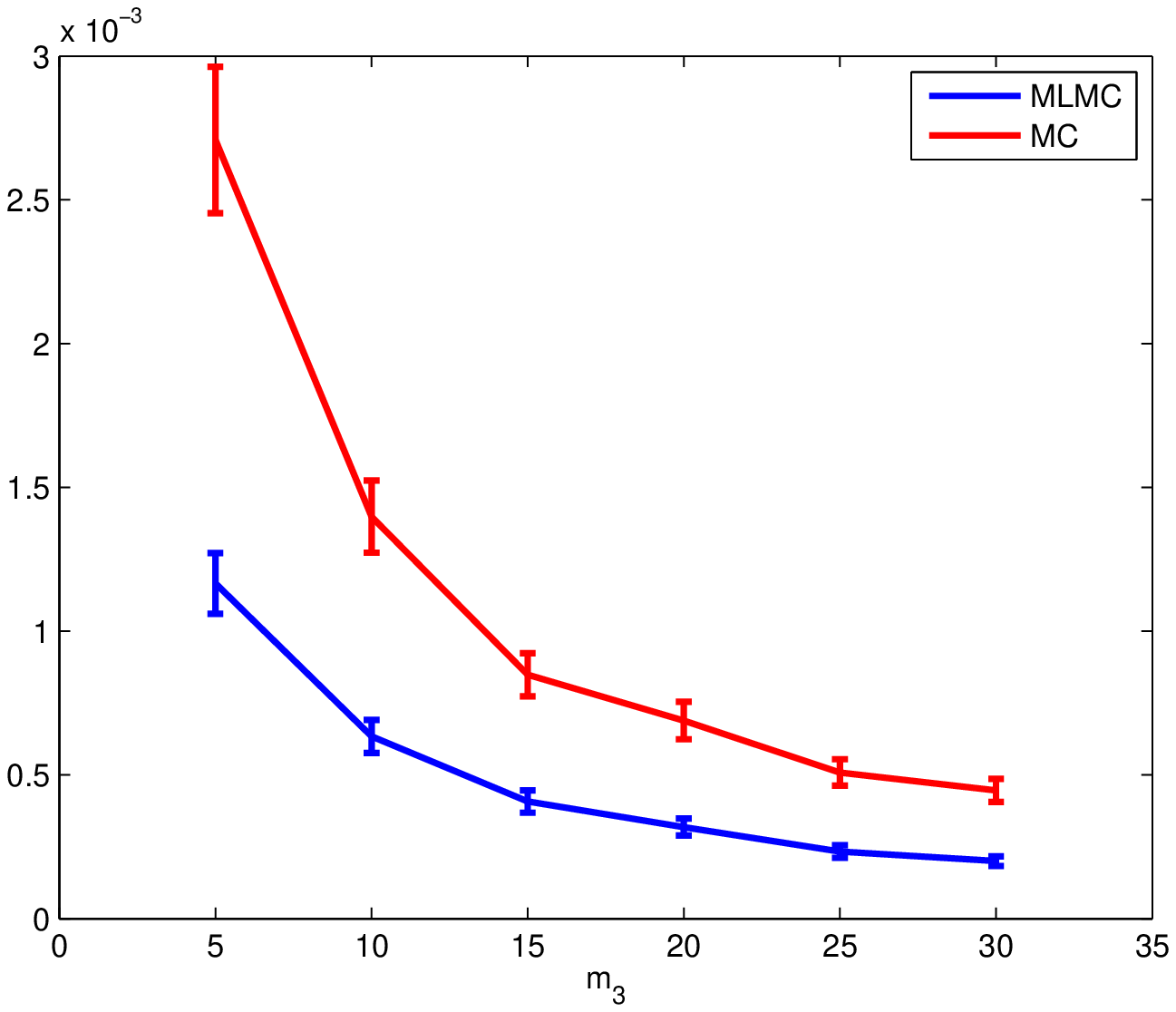}}}
\subfigure[Relative mean square errors on the two-point correlation of the effective coefficient.]
{\scalebox{.4}{ \includegraphics{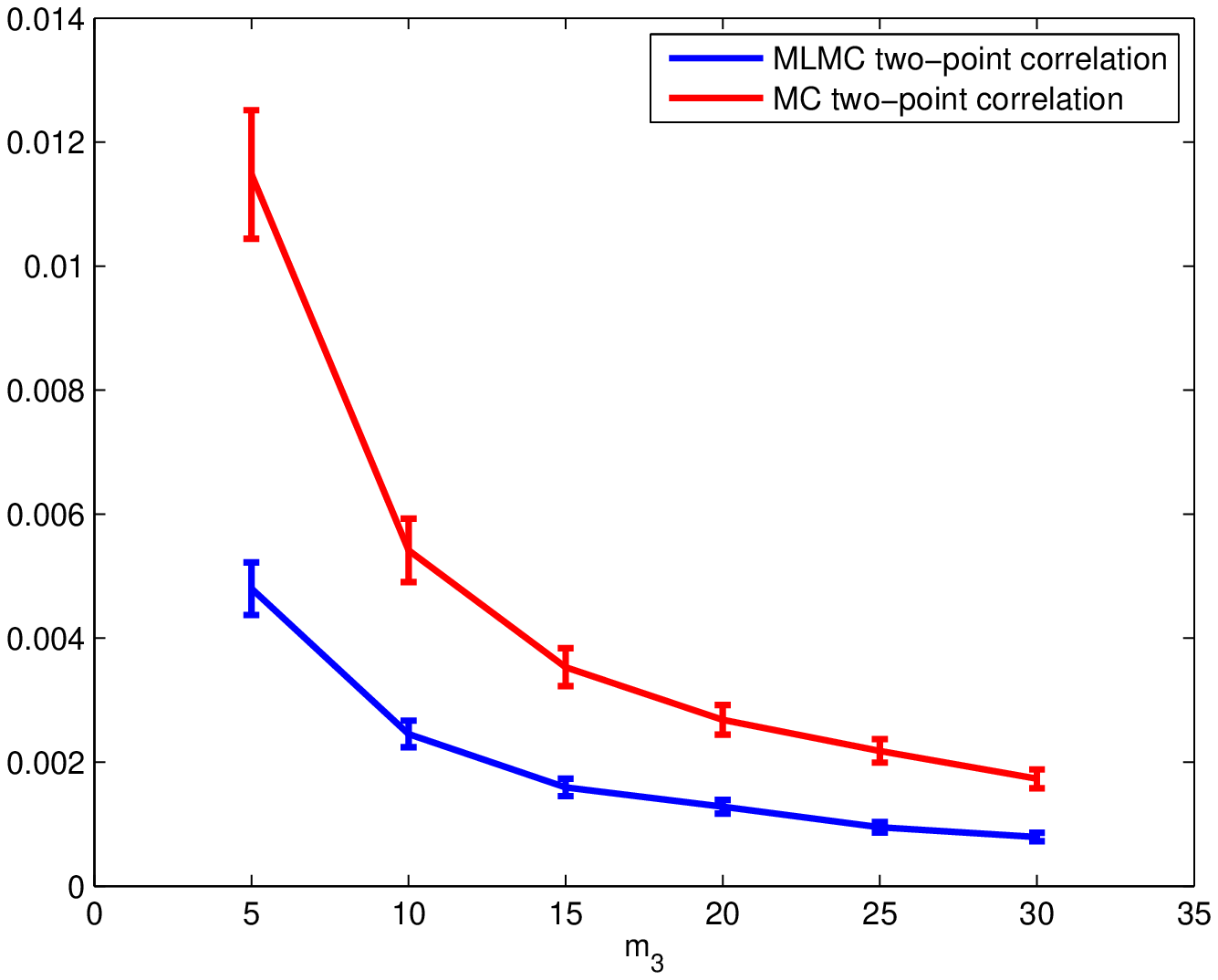}}}
\caption{Relative mean square errors with equated costs and $\mathfrak m
  =(16m_3, 4m_3, m_3)$, for the Example 1 (separable coefficient).}
\label{oned}
\end{figure}

For the MLMC approach, we use the RVEs $[a_l, b_l]=[0,0.5^{L+1-l}]$. For
this case, we expect that $\beta = 2$.
Following~\eqref{eq:m_l}, we hence take $\mathfrak m
=(4^{L-l}m_L,\cdots,4m_L,m_L)$ realizations. 
For comparison, we calculate the error of the standard MC approach on
the large RVE $[a_L,b_L] = [0, 0.5]$, with $\dis \widehat{m}_L =
\frac{\sum_{l=1}^L m_l \, b_l}{b_L}$ samples, so that both approaches
share the same cost. We are interested in comparing the relative mean
square errors 
$$
\left(e^{rel}_{MLMC}\right)^2(A^*_L) 
= 
\frac{\left(e_{MLMC}(A^*_L)\right)^2}{\left( \EE(A^*_L) \right)^2},
\quad
\left(e^{rel}_{MC}\right)^2(A^*_L)
= 
\frac{\left( e_{MC}(A^*_L) \right)^2}{\left( \EE(A^*_L) \right)^2},
$$
where $e_{MLMC}(A^*_L)$ and $e_{MC}(A^*_L)$ are defined by~\eqref{eq:e_mlmc}
and~\eqref{eq:e_mc}.
Since the errors depend on the set of chosen random numbers, we repeat
the computations $Nb=10000$ times
and calculate the corresponding confidence intervals for the errors:
$$
\left[ \text{mean}[(e^{rel})^2]-\frac{1.96 \
    \text{std}[(e^{rel})^2]}{\sqrt{Nb}},\text{mean}[(e^{rel})^2]+\frac{1.96 \ \text{std}[(e^{rel})^2]}{\sqrt{Nb}}\right].
$$
We take $L=3$, and show on Figure~\ref{oned} the relative mean square
errors on the expected value and the two-point correlation of the
effective coefficient. For both quantities, 
we observe that the MLMC approach yields errors 2.5 times smaller
than the MC approach at equal computational work.

\paragraph{Example 2 (separable stationary coefficient)}

We now consider an example where the effective coefficient does not
depend on $\omega'$, in the limit of infinitely large RVEs. We take $A$
with inverse given by
$$
A^{-1}(x, \omega, \omega')
=
\left(C +\sum_{i \in \ZZ} \chi_i (\omega') \, 1_{[i, i+1)}(x) \,
\sin^2({2 \pi x}) \right) \exp(\omega),
$$
where $\omega$ and $\chi_i$ are i.i.d. random variables, uniformly
distributed in $[0,1]$, $C=1$ and $1_{[i, i+1)}(x)$ denotes the
indicator function which is equal to $1$ for $x \in [i, i+1)$ and to zero
elsewhere. The apparent homogenized coefficient on the RVE $[a,b]$ (to
simplify, we choose $a$ and $b$ in $\ZZ$) is
\begin{eqnarray*}
A^*_{a,b}(\omega, \omega')
&=&\left( \frac{1}{b-a} \int_a^b A^{-1}(x, \omega, \omega') dx \right)^{-1}
\\
&=&\left( \frac{\exp(\omega)}{b-a}\left[ C(b-a)
+0.5 \sum_{i=a}^{b-1} \chi_i (\omega') \right]\right)^{-1}.
\end{eqnarray*}
In this case, the coefficient $A$ is stationary in the variables
$(x,\omega')$, hence the standard stochastic homogenization theory
holds: the exact effective coefficient is independent from $\omega'$
and reads
$$
A^*(\omega)
= 
\left[ \EE_{\omega'} \int_0^1 A^{-1}(x, \omega, \omega') dx \right]^{-1}
=
\left[ \exp(\omega) \left( C + 0.5 \EE(\chi) \right) \right]^{-1}.
$$
Remark that, as expected, $\dis \lim_{b-a \to \infty} A^*_{a,b}(\omega,
\omega') = A^*(\omega)$ almost surely in $\omega'$. In addition, the
Central Limit Theorem holds for this case, thus $\beta = 1$
in~\eqref{eq:hyp}. 

Following Section~\ref{sec:mlmc}, the theoretical reference value is
$\EE \left[ A^*_L(\omega,\omega') \right]$, where
$A^*_L(\omega,\omega')$ is the apparent homogenized coefficient on the
largest RVE. However, this theoretical reference value is not easy to
compute. We prefer to work with a different reference value, which is
analytically computable, and which is very close to $\EE \left[
  A^*_L(\omega,\omega') \right]$ when the RVE at level $L$ is large. 
In the sequel, we use as reference 
$$
A^*_{ref} := \EE \left[ A^*(\omega) \right]
=
\left( C + 0.5 \EE(\chi)\right)^{-1}
\EE \left[\exp(-\omega) \right]
=
\frac{1 - 1/e}{C + 0.25}.
$$
By construction, $\dis A^*_{ref} = \lim_{\eta_L \to \infty} 
\EE \left[ A^*_L(\omega,\omega') \right]$.

For the MLMC approach, we use the RVEs $[a_l, b_l]=[0,100 \times 2^{l-1}]$ with
$\mathfrak m =(2^{L-l}m_L,\cdots , 2m_L, m_L)$ realizations (recall that
$\beta = 1$ in this case, and hence this choice for $\mathfrak m$ agrees
with~\eqref{eq:m_l}). Note that the smallest RVE is again 
much larger than the characteristic length scale of the field $A$.
We compare this approach with a standard MC approach on the largest RVE 
$[a_L,b_L]=[0,100 \times 2^{L-1}]$ 
that uses $\dis \widehat{m}_L
=\frac{\sum_{l=1}^L m_l \, b_l}{b_L}$ samples (so that both approaches
share the same cost).

\begin{figure}[htp]\center
\subfigure[Relative mean square errors of the expected value of the effective coefficient ($L=3$).]
{\scalebox{.4}{ \includegraphics{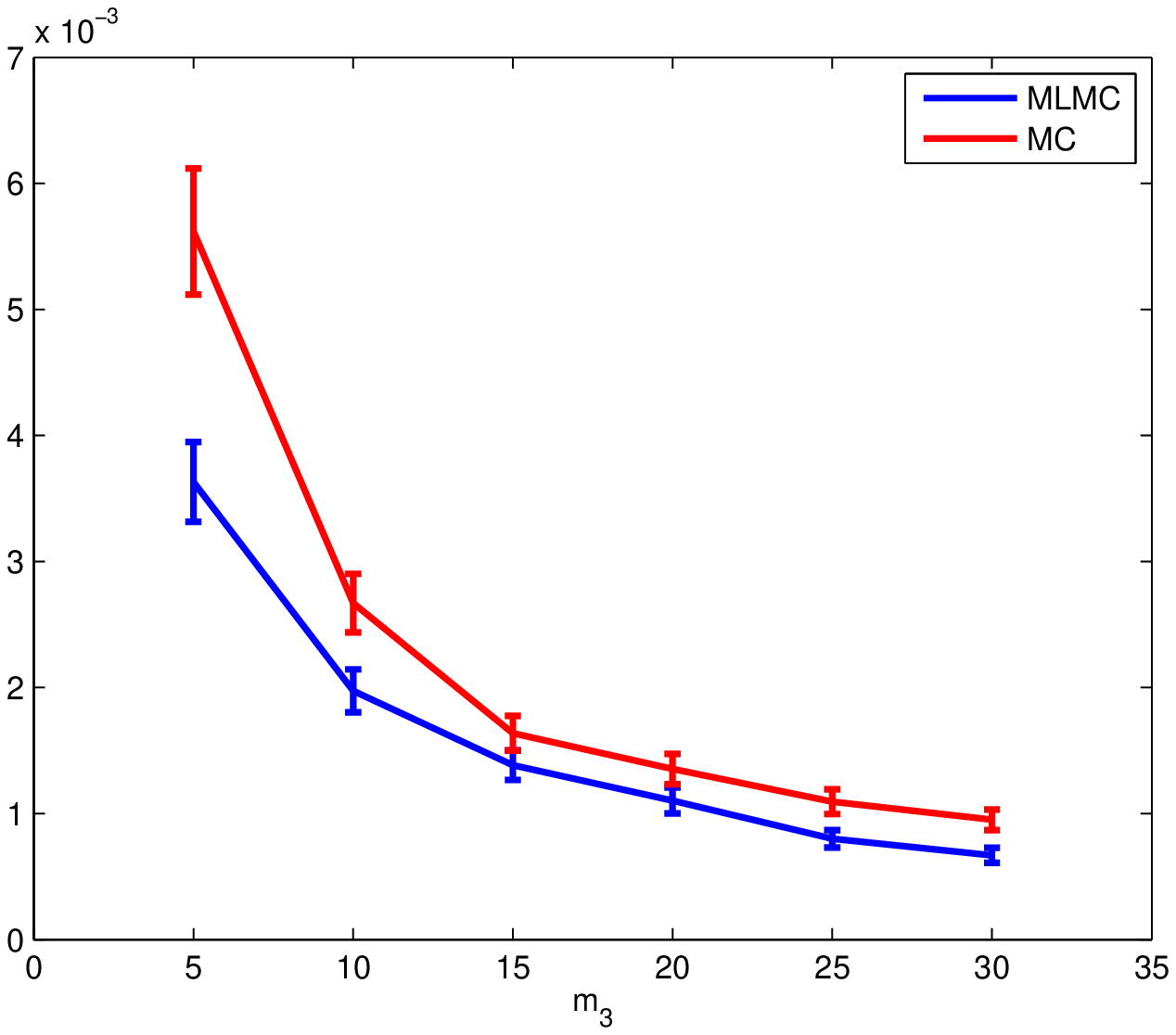}}}
\subfigure[Relative mean square errors of the two-point correlation of the effective coefficient ($L=3$).]
{\scalebox{.4}{ \includegraphics{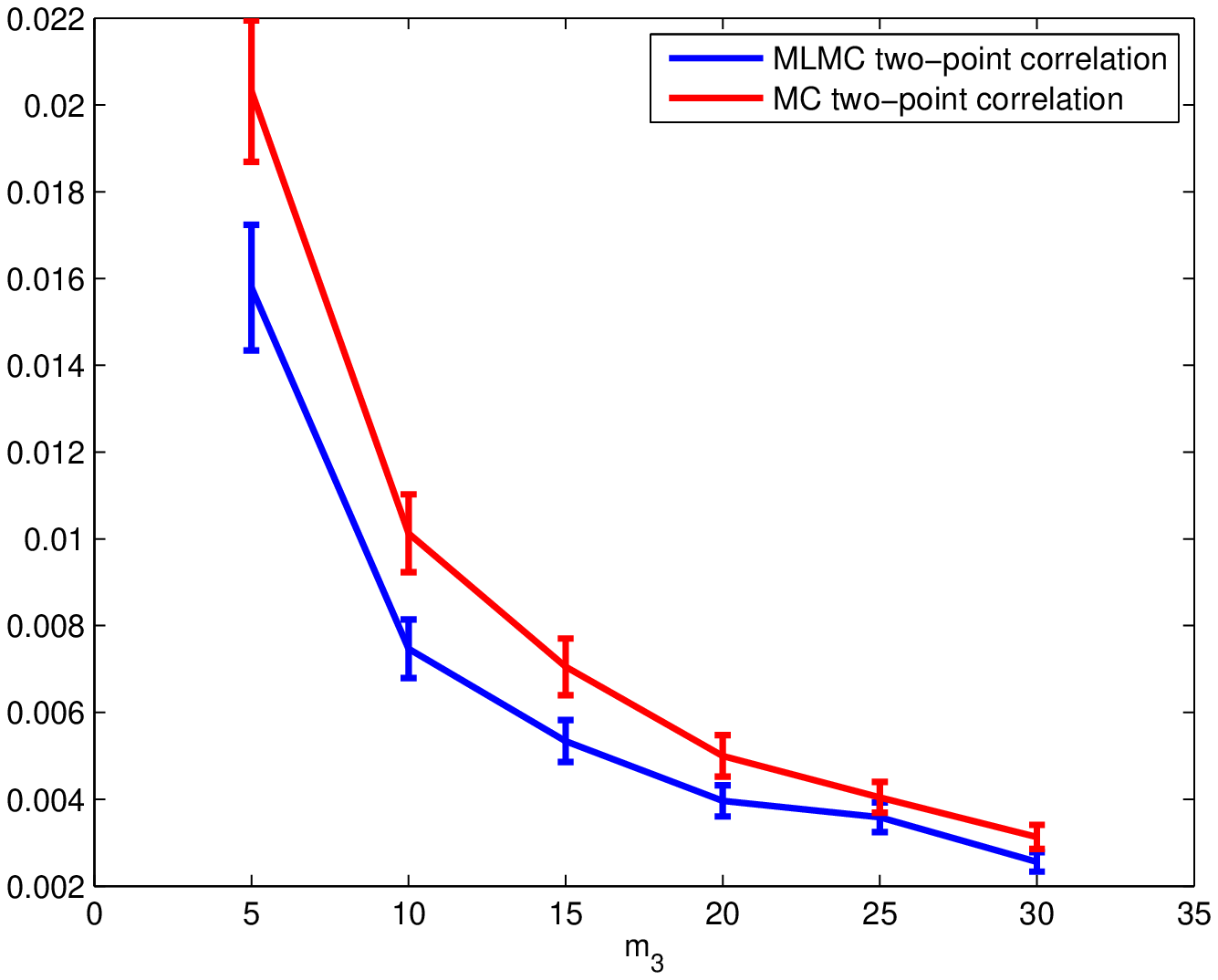}}}
\subfigure[Relative mean square errors of the expected value of the effective coefficient ($L=5$).]
{\scalebox{.4}{ \includegraphics{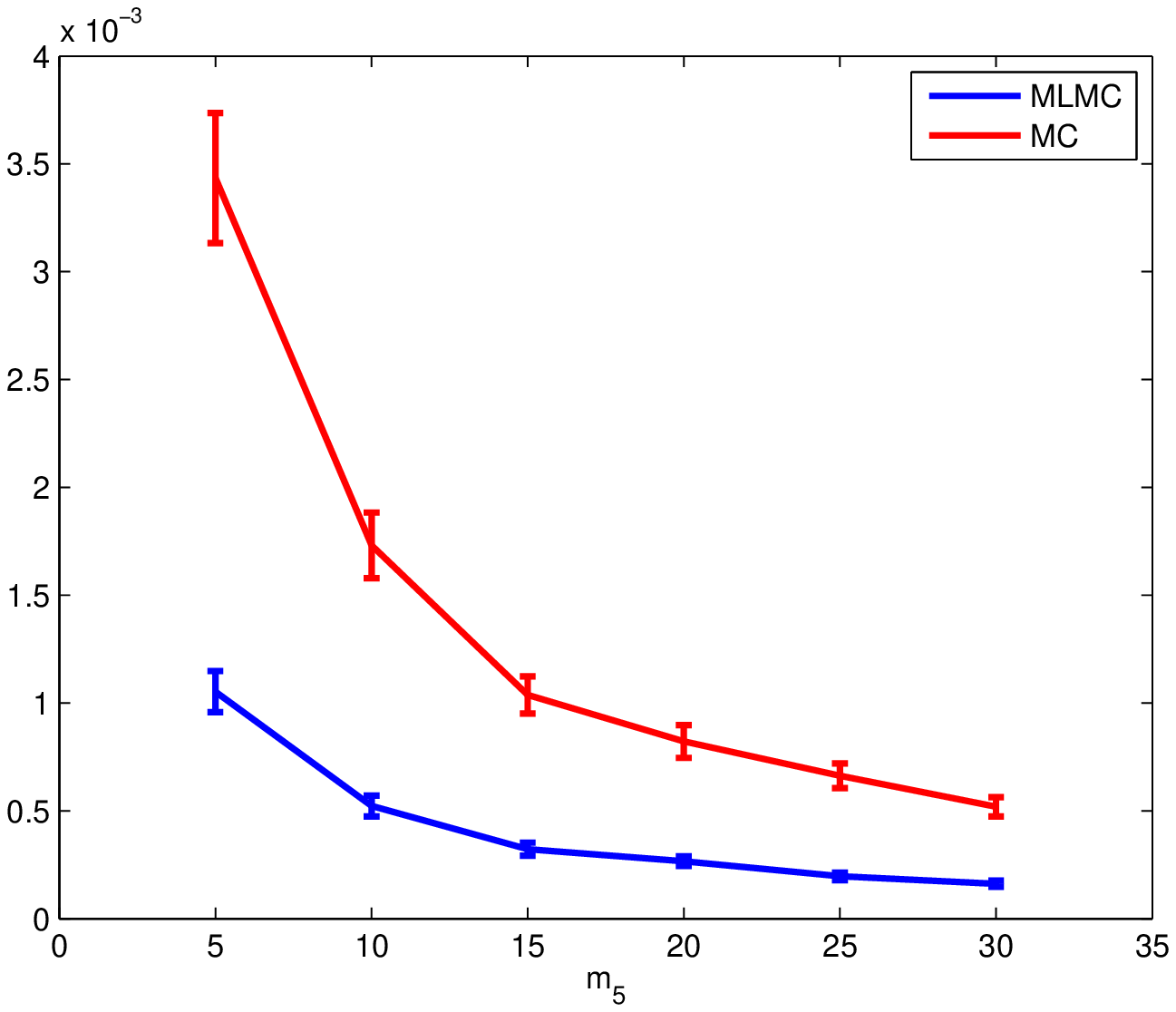}}}
\subfigure[Relative mean square errors of the two-point correlation of the effective coefficient ($L=5$).]
{\scalebox{.4}{ \includegraphics{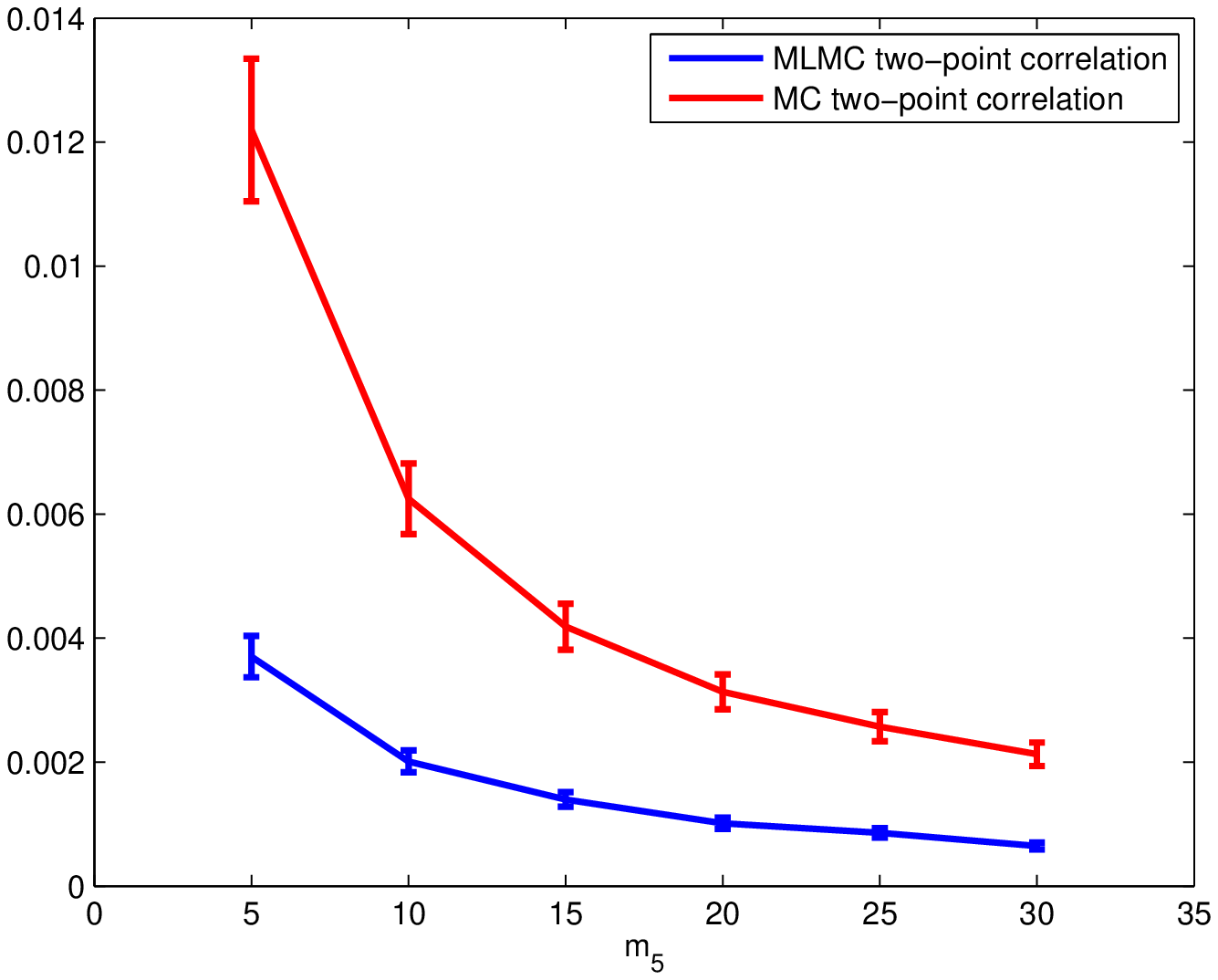}}}
\subfigure[Relative mean square errors of the expected value of the effective coefficient ($L=7$).]
{\scalebox{.4}{ \includegraphics{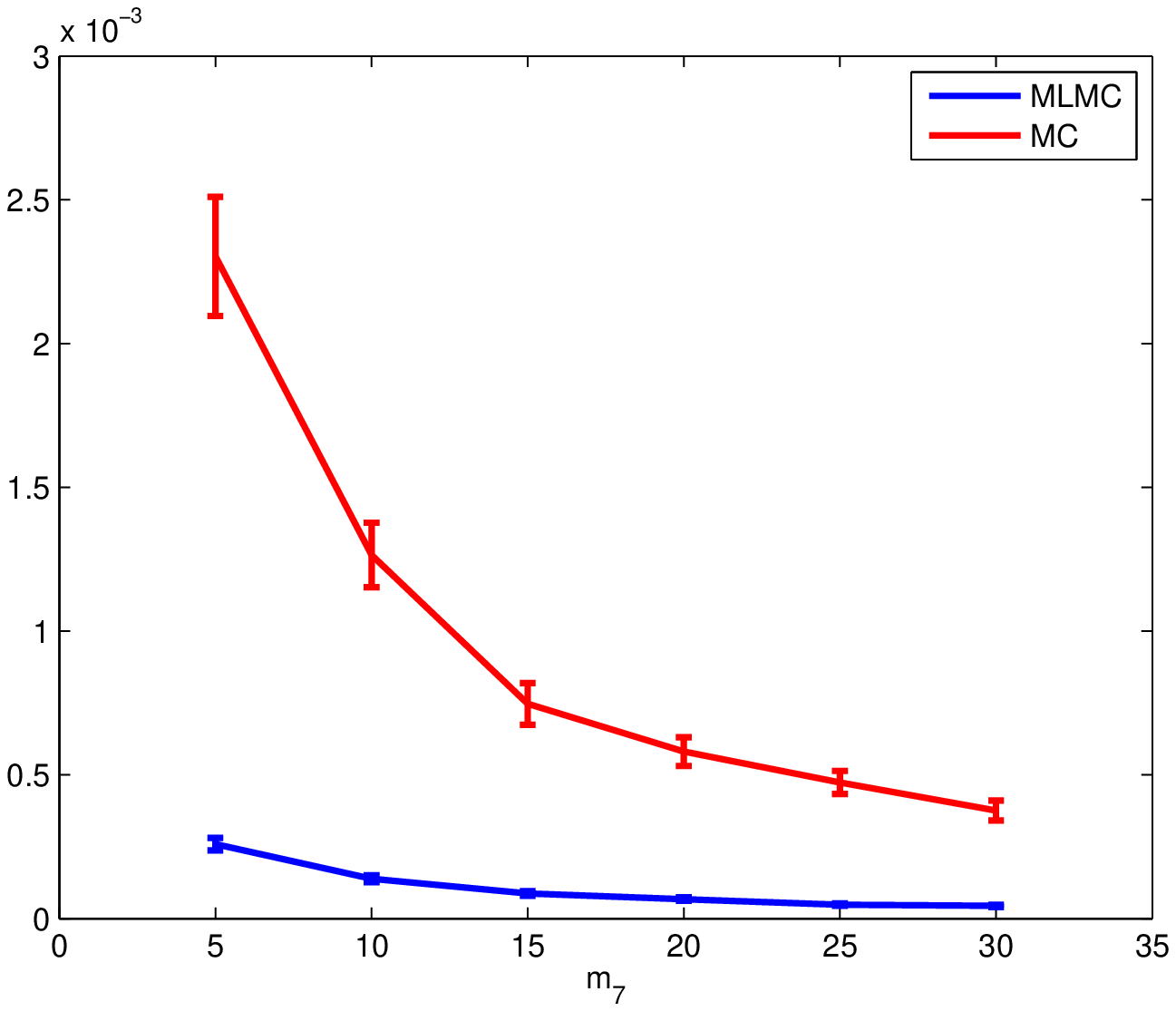}}}
\subfigure[Relative mean square errors of the two-point correlation of
the effective coefficient ($L=7$).]
{\scalebox{.4}{ \includegraphics{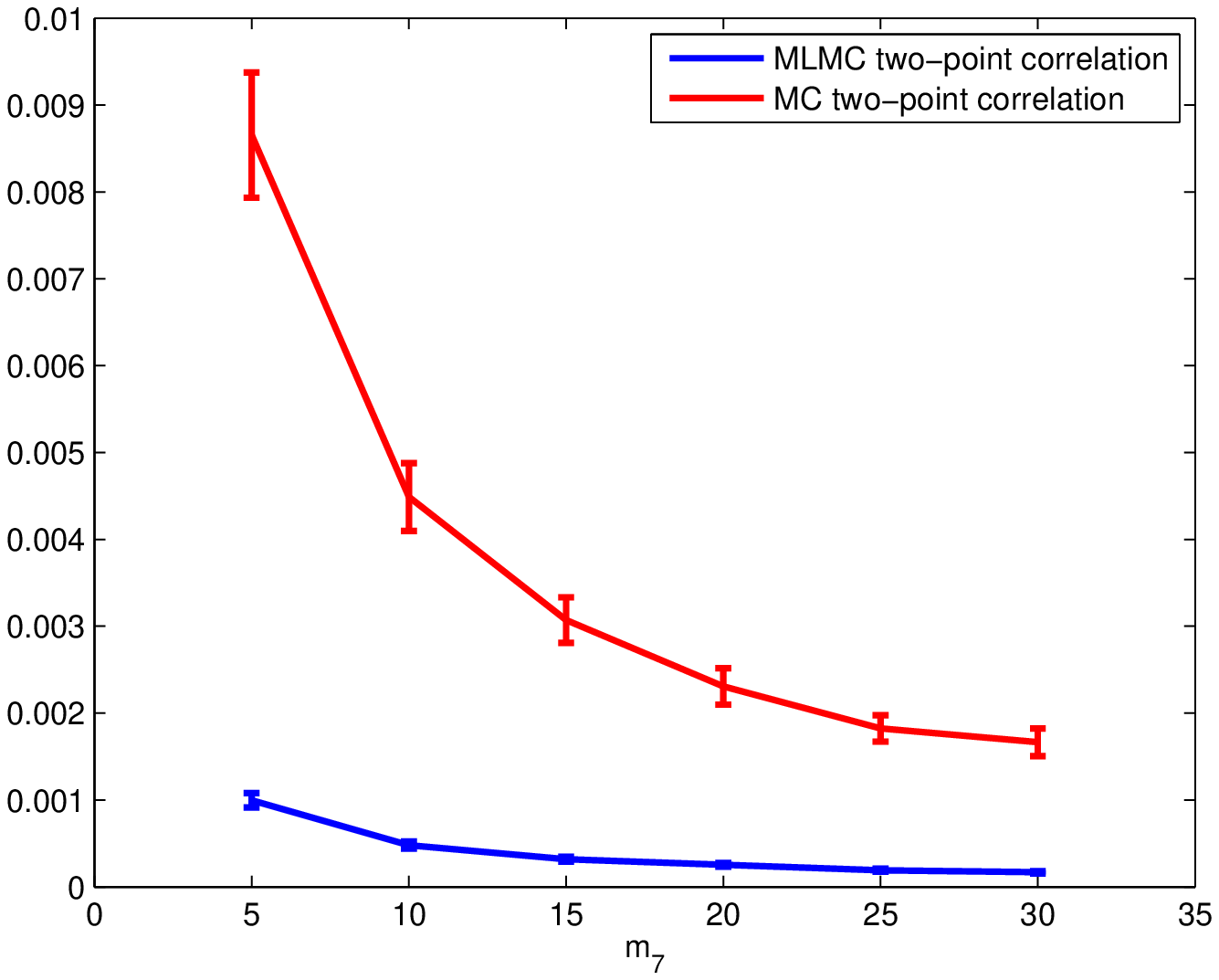}}}
\caption{Relative mean square errors with equated costs and 
$\mathfrak m =(2^{L-l}m_L,\cdots , 2m_L, m_L)$, for the Example
  2 (separable stationary coefficient).}
\label{onedstat}
\end{figure}

For this example, we have considered the choices $L=3$, 5 or 7.  
On Figure~\ref{onedstat}, we compare the relative mean square errors
$\left(e^{rel}_{MLMC}\right)^2$ and $\left(e^{rel}_{MC}\right)^2$ on the
expected value and the two-point correlation of the effective
coefficient (along with the corresponding confidence intervals obtained
from $Nb=10000$ different sets of random numbers).
We again observe that the MLMC approach is more accurate (for the same
amount of work), and that the gain in accuracy increases if we increase
the total number $L$ of levels (this observation is consistent with
Figure~\ref{rvework}). For $L=3$, the gain is equal to 1.5 for
both quantities, whereas it is equal to 3 for $L=5$ and to 8 when
$L=7$. 

\begin{figure}[htp]\center
{\scalebox{.4}{ \includegraphics{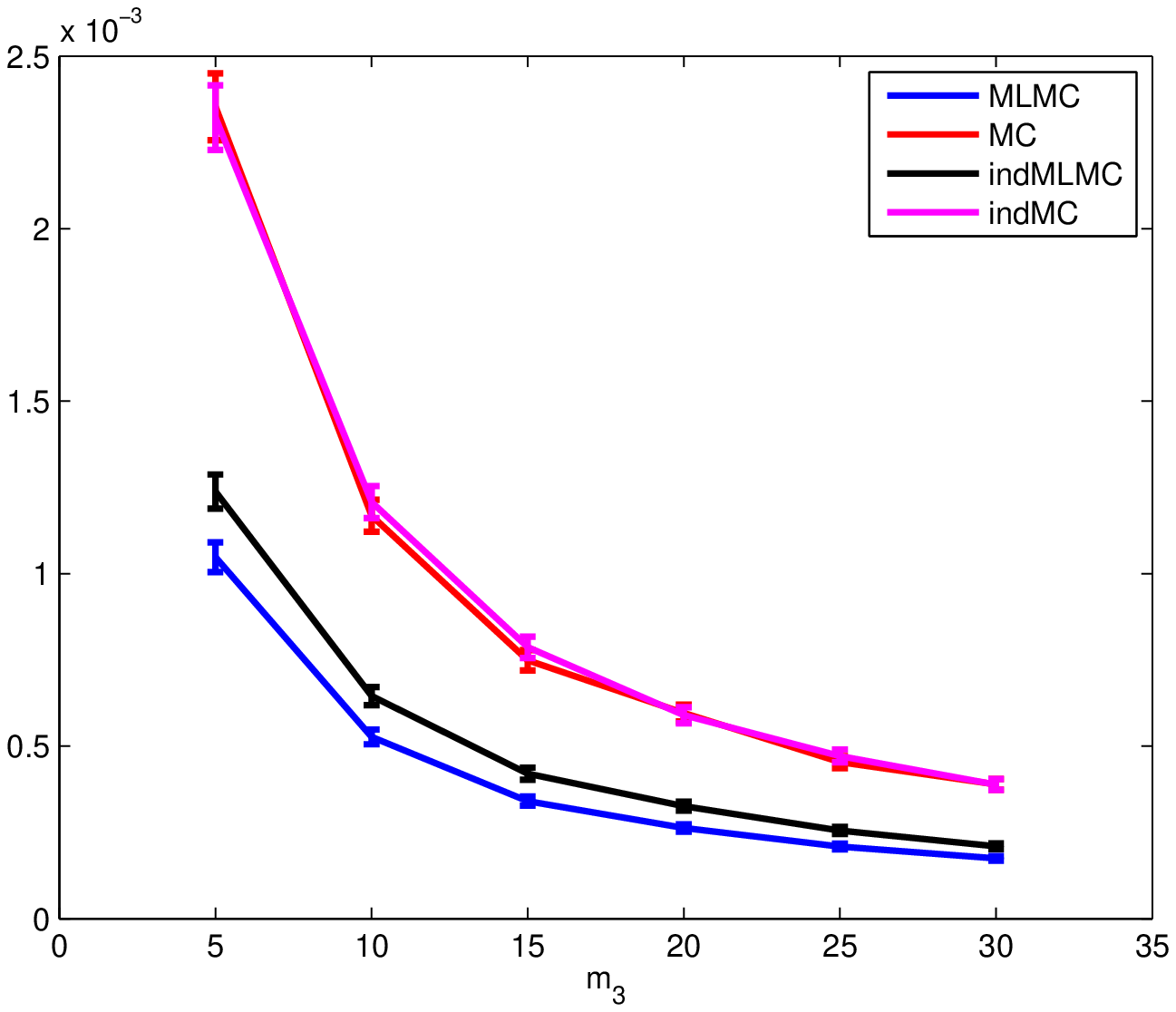}}}
\caption{Relative mean square errors for the Example 2 (separable
stationary coefficient). The MLMC results have been computed
following~\eqref{eq:mlmc_non_ind} with $(M_1,M_2,M_3) = (16 m_3,4m_3,
m_3)$. Results following the approach described in Remark~\ref{rem:iid},
labeled as 'MLMC ind', have been computed following~\eqref{eq:mlmc_ind}
with $(M_1,M_2,M_3) = (16 m_3,4m_3, m_3)$.
The MC results have been computed with $\widehat{m}_{MC}$ samples so
that the costs of 'MLMC', 'MLMC ind' and 'MC' are equal.
We work here with $M_{j-1} = 4 M_j$ (rather than $M_{j-1} = 2 M_j$ as in
Figure~\ref{onedstat}) to ensure that the number of samples per level
decreases. 
}
\label{fig:iid}
\end{figure}

\begin{bem}
\label{rem:iid}
On this example, we have also considered a MLMC approach where the
realizations of $X_l$ used in $E_{M_l}(X_l - X_{l-1})$ are independent
from the realizations of $X_l$ used in $E_{M_l}(X_{l+1} - X_l)$. More
precisely (assuming $L=3$ for the sake of simplicity), this approach
consists in approximating $\EE(X_{L=3})$ by
\begin{multline}
\label{eq:mlmc_ind}
E^{L=3}_{ind}(X_{L=3}) 
:= 
\frac{1}{M_3} \left[ \sum_{i=1}^{M_3} \left( X_3(\omega_i) -
  X_2(\omega_i) \right) \right]
\\
+
\frac{1}{M_2-M_3} \left[ \sum_{i=1+M_3}^{M_2} \left( X_2(\omega_i) -
  X_1(\omega_i) \right) \right]
+
\frac{1}{M_1-M_2} \sum_{i=1+M_2}^{M_1} X_1(\omega_i),
\end{multline}
rather than by
\begin{multline}
\label{eq:mlmc_non_ind}
E^{L=3}(X_{L=3}) 
= 
\frac{1}{M_3} \left[ \sum_{i=1}^{M_3} \left( X_3(\omega_i) -
  X_2(\omega_i) \right) \right]
\\
+
\frac{1}{M_2} \left[ \sum_{i=1}^{M_2} \left( X_2(\omega_i) -
  X_1(\omega_i) \right) \right]
+
\frac{1}{M_1} \sum_{i=1}^{M_1} X_1(\omega_i),
\end{multline}
as in~\eqref{eq:approx_mlmc}. We
compare on Figure~\ref{fig:iid} this method with a standard MC method,
where the number of samples has been chosen to again equate the costs. We
again observe that the MLMC approach~\eqref{eq:mlmc_ind} (with
independent samples) is more 
accurate than the MC approach. We also observe that, at equal cost, a
better accuracy is obtained when one uses~\eqref{eq:mlmc_non_ind} (with 
samples that are {\em not} necessarily independent) rather
than~\eqref{eq:mlmc_ind}.
\end{bem}

\paragraph{Example 3 (non separable coefficient)}

\begin{figure}[htp]
\center
\subfigure[Relative mean square errors of the expected value of the effective coefficient.]
{\scalebox{.4}{ \includegraphics{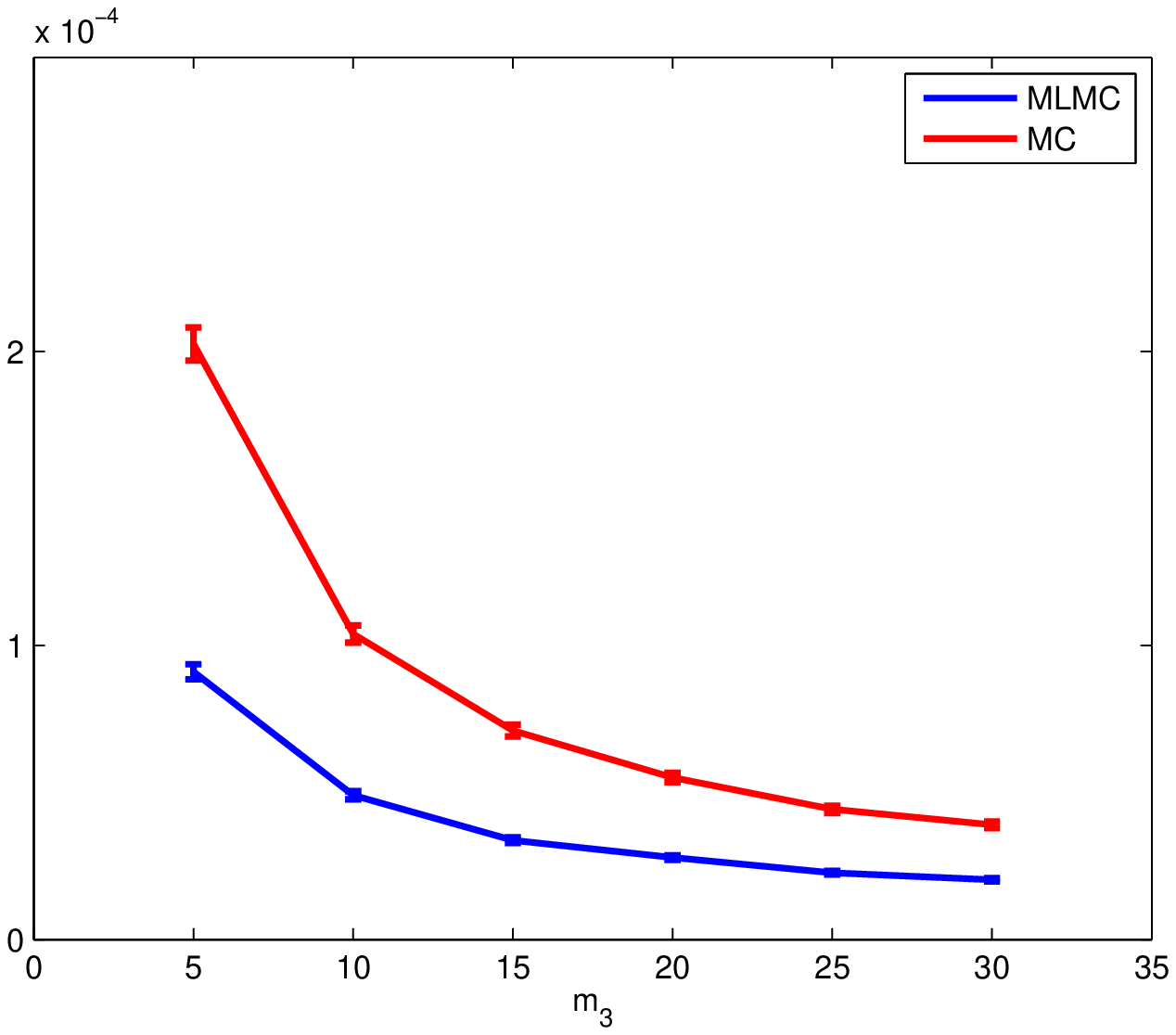}}}
\subfigure[Relative mean square errors of the two-point correlation of the effective coefficient.]
{\scalebox{.4}{ \includegraphics{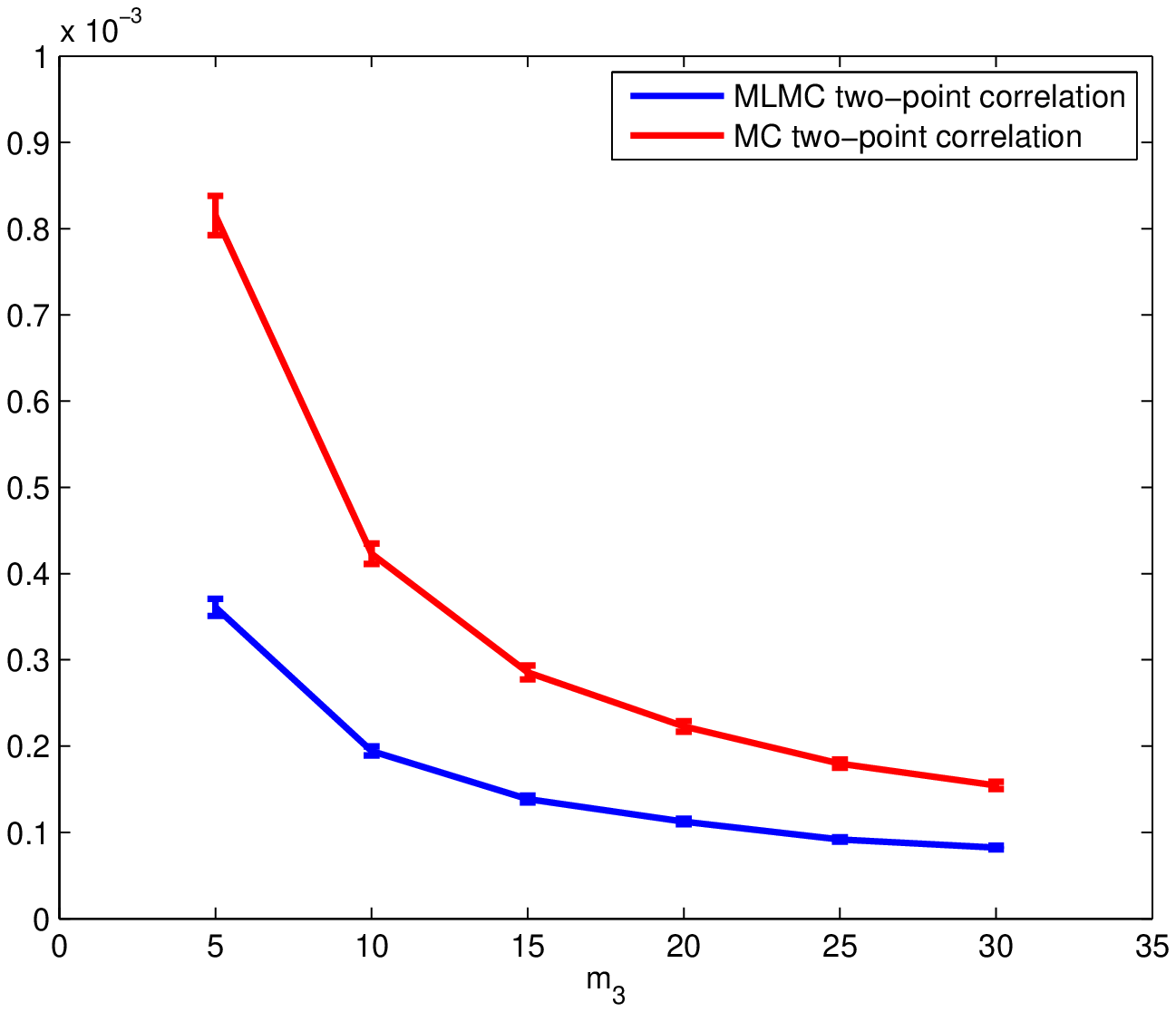}}}
\caption{Relative mean square errors with equated costs and $\mathfrak m
  =(16m_3, 4m_3, m_3)$, for the Example 3 (non separable coefficient).}
\label{onednonsep}
\end{figure}

We now consider the coefficient defined by its inverse as
$$
A^{-1}\left( \frac{x}{\epsilon}, \omega, \omega'\right) 
= 
C (1 + \omega)
+ 
\exp \left( \omega \omega' \sin\left(\frac{x}{\epsilon}\right) \right) 
\cos\left(\frac{x}{\epsilon}\right), 
$$
where $\omega$ and $\omega'$ are i.i.d. random variables uniformly
distributed in $[0.5, 1]$, $\dis \epsilon=\frac{0.5^L}{10}$ (the smallest RVE
is thus large compared to $\epsilon$) and $C=2e$ (which ensures that $A$ is
uniformly bounded away from 0). This example is more challenging than
the two previous ones as it is not separable. 
The apparent effective coefficient on $[a,b]$ is
\begin{eqnarray*}
& & A^*_{a,b}(\omega, \omega')\\
&=&\left( \frac{1}{b-a} \int_a^b 
A^{-1}\left(\frac{x}{\epsilon}, \omega, \omega' \right) \right)^{-1}
\\
&=&\left( \frac{1}{b-a}
   \left[ C(1+\omega) (b-a)
      + \frac{\epsilon}{\omega \omega'} \left( 
\exp\left( \omega \omega' \sin\left(\frac{b}{\epsilon}\right) \right) 
 -\exp\left( \omega \omega' \sin\left(\frac{a}{\epsilon}\right) \right) 
\right) \right]\right)^{-1}.
\end{eqnarray*}
As for the previous example, we use the practical reference value
$$
A^*_{ref} := \lim_{\eta \to \infty} \EE \left[ 
A^*_\eta(\omega, \omega') \right]
=
\frac{2 \ln 4/3}{C}.
$$
As for Example 1, we expect in this case that $\beta = 2$ and use
the RVEs $[a_l,b_l]=[0,0.5^{L+1-l}]$ with 
$\mathfrak m =(4^{L-l}m_L,\cdots,4m_L,m_L)$ realizations for the MLMC
approach, and compare its accuracy (at equal cost) with MC results on
the RVE $[a_L,b_L] = [0, 0.5]$.
Choosing $L=3$, we show on Figure~\ref{onednonsep} the relative mean
square errors $\left(e^{rel}_{MLMC}\right)^2$ and
$\left(e^{rel}_{MC}\right)^2$ on the expected value and the two-point
correlation of the effective coefficient (confidence intervals have
again been obtained from $Nb=10000$ different sets of random numbers).
Again, the MLMC approach yields an accuracy gain (here of the order of
2) over the MC approach, for both quantities. 

\subsubsection{Two dimensional examples}
\label{sec:coeff_2D}

We have seen in the previous section that the MLMC approach is efficient
in the one dimensional case. We turn here to two
dimensional test cases. 

\paragraph{Example 1 (separable coefficient)}

\begin{figure}[htp]
\center
\scalebox{.4}{ \includegraphics{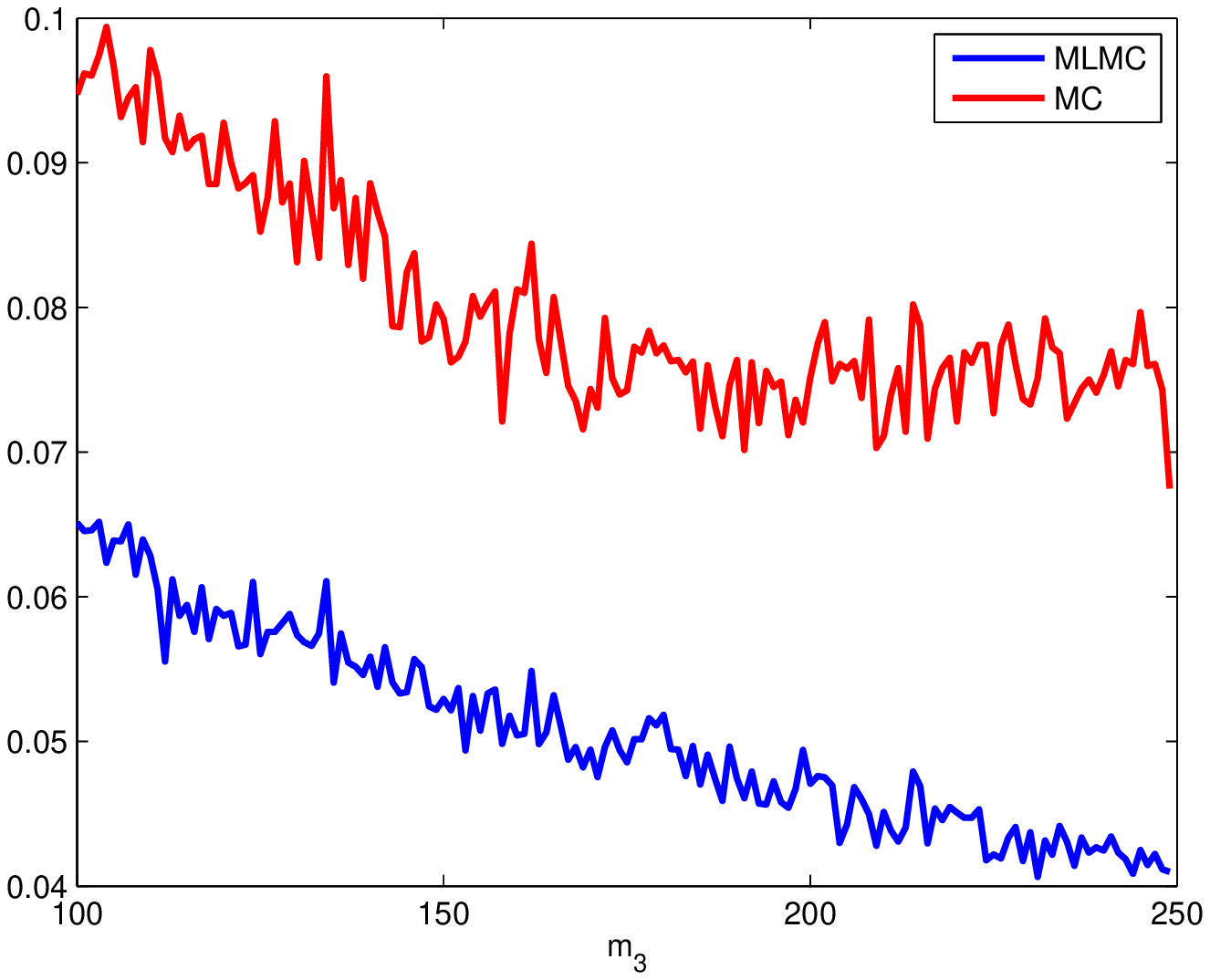}}
\caption{Errors $e_{MC}(\left[A^*_L\right]_{11})$ and
  $e_{MLMC}(\left[A^*_L\right]_{11})$ for $\mathfrak m=(4 m_3, 2 m_3,
  m_3)$, for the Example 1 (separable coefficient).}
\label{coeff}
\end{figure}

We first study the case when there is a separation in the randomness at
the macroscopic level and the microscopic level. We set
$$
A\left(\omega, \frac{x}{\epsilon}, \omega'\right)=
\widetilde{A}(\omega) \, B\left(\frac{x}{\epsilon}, \omega'\right),
$$
where $\widetilde{A}$ and $B$ are both scalar valued. The random field
$B$ has expected value $\EE(B)=10$ and a Gaussian covariance function:
$$
\text{cov}(x,x') = \Cov \left[ B \left( \frac{x}{\epsilon}, \cdot \right),
B \left( \frac{x'}{\epsilon}, \cdot \right) \right]
=
\sigma^2 \exp \left( -\frac{ | x-x' |^2}{\epsilon^2 \tau_0^2} \right),
$$ 
with $\sigma=\sqrt{2}$, $\tau_0 = \sqrt{2}$ and $\tau = \epsilon \tau_0 =
0.04$. We generate samples of the coefficient with the Karhunen-Lo\`eve
expansion. 
We take $\widetilde{A}(\omega)= \exp( \omega)$,
where $\omega$ is distributed according to the Gaussian law $N(0,1)$.
The effective matrix is $A^*(\omega,\omega')=\widetilde{A}(\omega) \,
B^*(\omega')$. We only define levels $l$ to approximate the expectation
of $B^*$. Thus, at each level $l$, we define
$A^*_l(\omega,\omega')=\widetilde{A}(\omega) \, B_l^*(\omega')$. 
Using $m_l$ independent samples at the microscopic level 
$\left\{ \omega'_j \right\}_{1 \leq j \leq m_l}$ and 
$n \times m_l$ independent samples at the macroscopic level 
$\left\{ \omega^i_j \right\}_{1 \leq j \leq m_l, \, 1 \leq i \leq n}$, we
define, for any $1 \leq i \leq n$,
$$
E_{m_l}(A^*_l)(\omega^i)
:= 
\frac{1}{m_l} \sum_{j=1}^{m_l}
A^*_l(\omega^i_j,\omega'_j)
= 
\frac{1}{m_l} \sum_{j=1}^{m_l}
\widetilde{A}(\omega^i_j) \, B_l^*(\omega'_j). 
$$
To approximate expectations at the microscopic level, we use the MLMC
approach, and introduce, for any $1 \leq i \leq n$,
\begin{eqnarray*}
E^L(A^*_L)(\omega^i)
&:=& 
\sum_{l=1}^L E_{m_l}(A^*_l-A^*_{l-1})(\omega^i)
\\
&=& 
\sum_{l=1}^L \frac{1}{m_l} \sum_{j=1}^{m_l} \widetilde{A}(\omega^i_j) 
\, ( B_l^*(\omega'_j)-B_{l-1}^*(\omega'_{j})).
\end{eqnarray*}
Expectations at the macroscopic level are approximated using a
standard MC approach on the macroscopic random variable $\omega$. The
reference quantity we are after is
$$
A^*_{ref} = \frac{1}{L} \sum_{l=1}^L \EE(A^*_l),
$$
which is in practice approximated by
$$
A^*_{ref} = \frac{1}{L} \sum_{l=1}^L \frac{1}{n} \sum_{i=1}^n
\frac{1}{m_l^{ref}} \sum_{j=1}^{m_l^{ref}} 
\widetilde{A}(\omega^i_j) \, B_l^*(\omega'_j). 
$$
In this case, the errors read
\begin{eqnarray*}
e_{MLMC}(A^*_L) 
&=& 
\sqrt{\frac{1}{n} \sum_{i=1}^{n} \left[
A^*_{ref} - 
E^L\left(A^*_L\right)\left(\omega^i\right) \right]^2},
\\ 
e_{MC}(A^*_L) 
&=& 
\sqrt{\frac{1}{n} \sum_{i=1}^{n} \left[
A^*_{ref} - 
E_{\widehat{m}_L}(A^*_L)(\omega^i) \right]^2},
\end{eqnarray*}
where $\widehat{m}_L$ is the number of samples used in the MC approach. 
As mentioned above, we equate the computational work of the MLMC
approach, which is
$\dis \sum_{l=1}^L m_l \left(\frac{\eta_l}{\epsilon}\right)^2$, with
that of the MC approach, which is $\dis \widehat{m}_L
\left(\frac{\eta_L}{\epsilon}\right)^2$. This leads to taking
$$
\widehat{m}_L = 
\frac{\sum_{l=1}^L m_l \left(\eta_l/\epsilon \right)^2}{\left(\eta_L/\epsilon\right)^{2}}.
$$
We next compare the errors. We choose to work with $L=3$ levels, 
$n=500$ and, to compute the reference value $A^*_{ref}$, we used 
${\mathfrak m}^{ref}=(m_1^{ref},m_2^{ref},m_3^{ref}) = (2000,\,1000,
\,300)$. We also adopt the parameters of Table~\ref{data}.
On Figure~\ref{coeff}, we show the errors on the first entry of the
effective matrix, $e_{MC}(\left[A^*_L\right]_{11})$ and
$e_{MLMC}(\left[A^*_L\right]_{11})$, for $\mathfrak m=(4 \, m_3, 2
\, m_3, \, m_3)$ (this choice is consistent with the value $\beta = 1$
in~\eqref{eq:hyp}; in turn, this assumption for $\beta$ is consistent
with our empirical estimation detailed in Section~\ref{sec:estim_beta}).
We observe from these simulations that the
MLMC approach provides (roughly twice as) smaller errors than the
standard MC approach for the same amount of computational work. Similar
conclusions hold for the other entries of the effective matrix. 

\paragraph{Example 2 (non-separable coefficient)}

We now consider a more difficult case, where there is no separation
between uncertainties at the macro- and the microscopic levels. In
general, such cases are difficult to handle, since having a sufficiently
large number of samples to appropriately reduce the statistical noise is
very expensive. We consider 
below a specific example for $A$ such that we can solve the local
problems (and thus compute the effective coefficient) analytically, due
to the specific choice of boundary conditions in the local problem. Note
that, in the limit of 
infinitely large RVEs, the effective coefficient does not depend on the
precise choice of the boundary conditions set on the local problems
(see~\cite{BP04}).

We consider the scalar coefficient ($x=(x_1,x_2)$)
$$
A\left(x,\omega,\frac{x}{\epsilon}, \omega'\right)
=
A_1\left(x_1,\omega,\frac{x_1}{\epsilon},\omega'\right) \ 
A_2\left(x_2,\omega, \frac{x_2}{\epsilon}, \omega'\right),
$$
and write the local problems with Dirichlet and no-flow boundary
conditions: 
\begin{eqnarray*}
-\text{div}\left( A\left(x,\omega,\frac{x}{\epsilon},\omega'\right)
  \nabla \chi_i\right) 
&=& 
0 \text{ in $Y_\eta = (0,\eta)^d$},
\\
\chi_i(x,\omega,\omega') &=& x_i \text{ on $\partial Y^D_\eta$},
\\
n\cdot \nabla \chi_i &=& 0 \text{ on $\partial Y_\eta \setminus \partial Y^D_\eta$},
\end{eqnarray*}
with $\partial Y^D_\eta = 
\lbrace x \in \partial Y_\eta \, | \, x_i=0 \text{ or } x_i=\eta\rbrace$.
With these choices, the local problem reduces to a one-dimensional
problem in the direction $x_i$ for the function $\chi_i$ that only
depends on $x_i$. For the first entry of the upscaled coefficient, we get
\begin{equation}
\label{eq:decoupl_dim}
A^*_{11}(\omega, \omega')
=
\left( \frac{1}{\eta} \int_0^\eta A_1^{-1}
\left(x_1,\omega,\frac{x_1}{\epsilon},\omega'\right)
\, dx_1 \right)^{-1} 
\frac{1}{\eta} \int_0^\eta
A_2\left(x_2,\omega,\frac{x_2}{\epsilon},\omega'\right) \, dx_2.
\end{equation}
In our example we choose
\begin{eqnarray*}
A_1^{-1}\left( \omega, \frac{x_1}{\epsilon}, \omega'\right)
&=& 
C (1 + \omega)+ 
\exp\left(  \omega \omega' \sin\left(\frac{x_1}{\epsilon}\right) \right) 
\cos\left(\frac{x_1}{\epsilon}\right), 
\\
A_2\left(x_2,\omega, \frac{x_2}{\epsilon},\omega'\right)
&=& 
C (1 + \exp(5\omega)) \, x_2 
\\
&& + 
(1+x_2) \exp\left( (1+x_2) \omega \omega'
  \sin\left(\frac{x_2}{\epsilon}\right) \right) \cos
\left(\frac{x_2}{\epsilon} \right), 
\end{eqnarray*}
where $\omega$ and $\omega'$ are i.i.d. random variables uniformly
distributed in $[0.5,1]$ and $C=2e$. In this case, we see
that~\eqref{eq:hyp} holds with $\beta = 2$.
To ensure scale
separation even for the smallest RVE, we take 
$\dis \epsilon= \frac{\eta_1}{10}$.

To define the reference value of the effective coefficient, we run a MC
approach on the RVE $[0,0.5]^2$ with $\widetilde{m}= 400000$
realizations. It is possible to compute such a large number of samples
in this two-dimensional test case thanks to the specific
analytical expression~\eqref{eq:decoupl_dim}. 

\begin{figure}[htp]
\center
\subfigure[Relative mean square errors of the expected value of the effective coefficient.]{
\scalebox{.4}{ \includegraphics{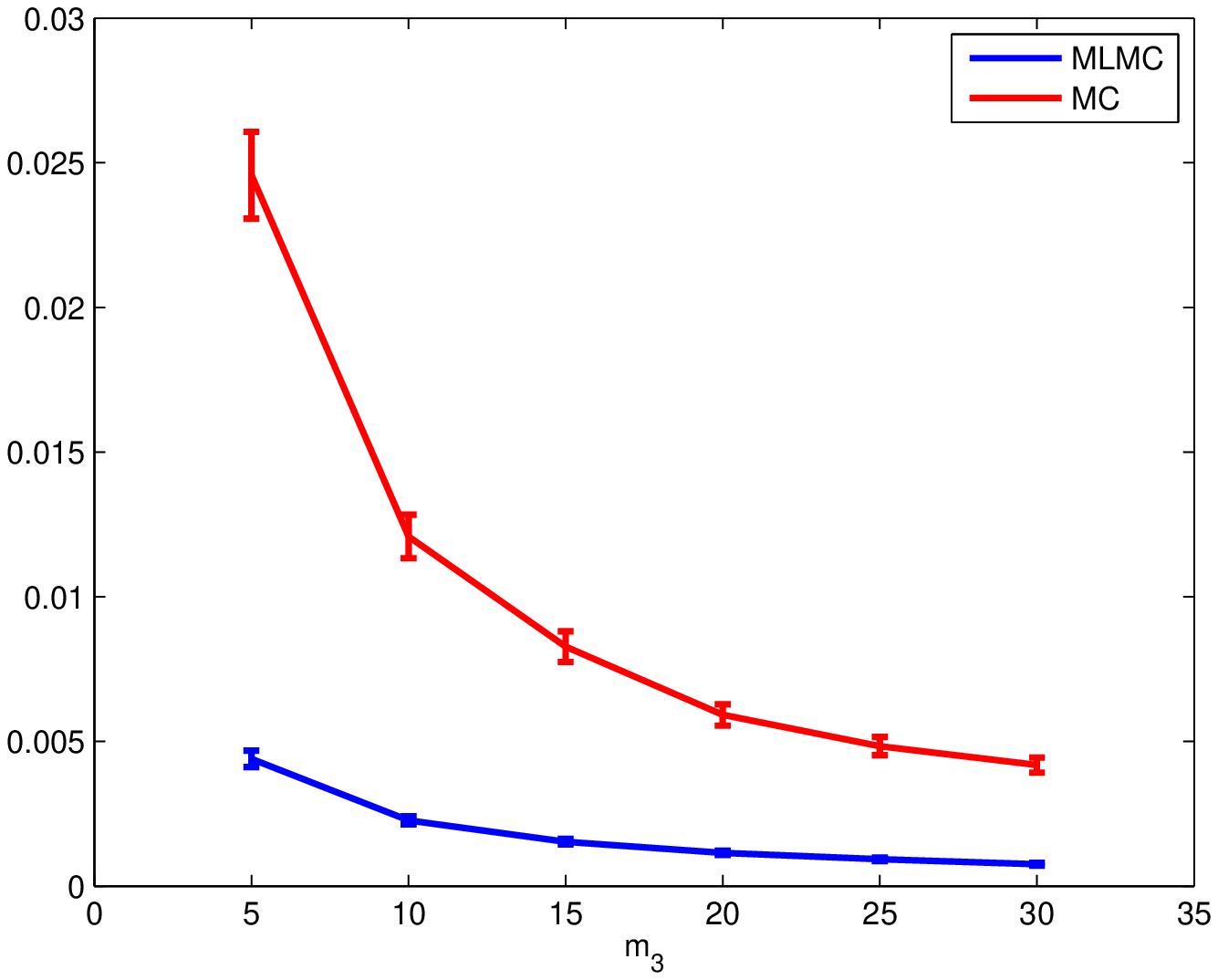}}}
\subfigure[Relative mean square errors of the two-point correlation of the effective coefficient.]{
\scalebox{.4}{ \includegraphics{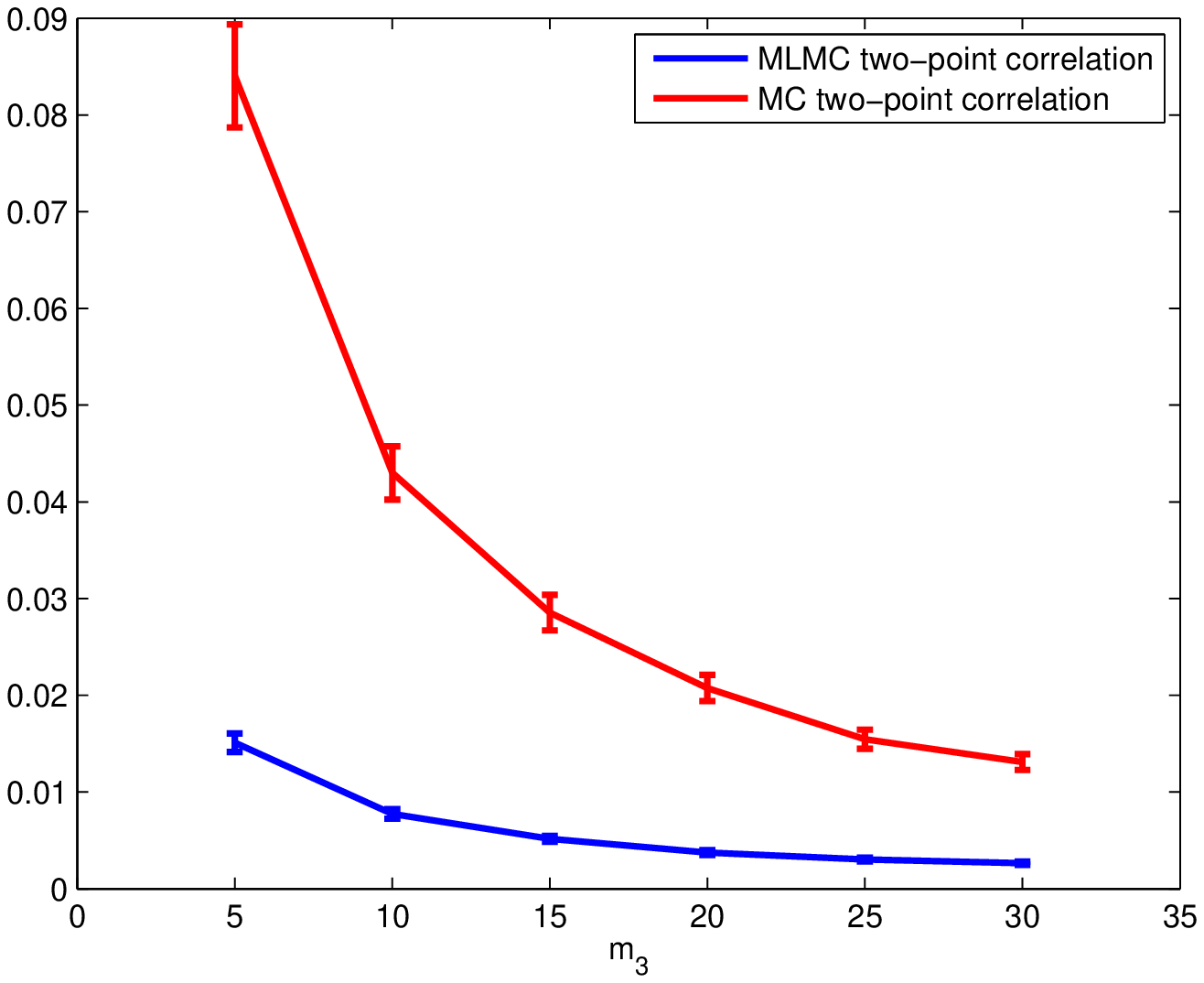}}}
\caption{Relative mean square errors with equated costs and $\mathfrak m
  =(16m_3, 4m_3, m_3)$, for the Example 2 (non-separable coefficient).}
\label{2dnonsep}
\end{figure}

The MLMC approach is run with $L=3$ different levels, and $\mathfrak m=
(16 m_3, 4 m_3, m_3)$ realizations at each level (a choice which is
consistent with~\eqref{eq:m_l} and the fact that $\beta=2$). To determine a
confidence interval, we repeat the overall procedure with $Nb=2000$
different sets of realizations. We compare on Figure~\ref{2dnonsep} the
accuracies of the MC and MLMC approaches at equal computational
cost. Again, the MLMC approach is more accurate, here by a factor roughly
equal to 5.

\subsection{Computation of the homogenized  solution}
\label{sec:num_sol}

\subsubsection{One dimensional example}

As in Section~\ref{sec:num_coeff}, we start with the one dimensional
situation where we know the 
reference solution exactly. To make the computations even simpler, we
assume that, at the coarse-scale, the problem is subjected to
homogeneous Neumann boundary conditions. The coarse problem thus reads 
\begin{equation}
\label{eq:edp_1D}
\frac{d}{dx}\left(A^*(x, \omega, \omega')\frac{d u^*}{dx} \right)= f(x),
\quad
(u^*)'(0) = (u^*)'(1) = u^*(0) = 0,
\end{equation}
where the right-hand side satisfies $\dis \int_0^1 f=0$. The exact
solution is
$$
u^*(x,\omega,\omega')= \int_0^x \left(A^*(t,\omega,\omega')\right)^{-1}
F(t) \, dt,
\quad F(t) = \int_0^t f(z) \,dz.
$$
Let $x_i$ denote the vertices of the grid, $0 \leq i \leq N$. 
The numerical approximation of $u^*$ is a piecewise constant function,
equal, on the interval $(x_{i-1},x_i)$, to
$$
u^*_i = \sum_{j=1}^i (A^*(x_j,\omega,\omega'))^{-1} \int_{x_{j-1}}^{x_j} \, F(x) \,dx. 
$$
In the spirit of the Example 3 in Section~\ref{sec:coeff_1D}, we assume
that the apparent homogenized coefficient, obtained by solving the local
RVE problem on $[a,b]$, reads
\begin{multline*}
(A^*(x, \omega, \omega'))^{-1}
=
C(1+\exp(5\omega))x 
\\
+ \frac{1}{b-a} \frac{\epsilon}{\omega \omega'}
\left[
           \exp\left( (1+x) \omega \omega' \sin\left(\frac{b}{\epsilon}\right) \right) 
           -\exp\left( (1+x) \omega \omega' \sin\left(\frac{a}{\epsilon}\right)
           \right) 
\right]
\end{multline*}
where $\omega$ and $\omega'$ are i.i.d. random variables uniformly
distributed in $[0.5,1]$ and $C=2e$. We take $f(x)=e^x-(e-1)$.

The reference quantity is the expectation of the solution
to~\eqref{eq:edp_1D}, computed with the coefficient $A_\infty^*$ obtained by
considering an infinitely large RVE:
$$
(A_\infty^*(x, \omega, \omega'))^{-1}
=
C(1+\exp(5\omega))x. 
$$
This reference quantity reads
$$
\EE\left(u^*_\infty\right)
=
{\cal C} \left( \exp(x)x-\exp(x)-(e-1)\frac{x^3}{3}-\frac{x^2}{2} +1 \right)
$$
where $\dis {\cal C} = C \left( 1 + 2 \int_{1/2}^1 e^{5 \omega} d\omega
\right)$.

To compute an approximation of $\EE\left(u^*_\infty\right)$, we use the
MLMC approach with $L=3$ levels.
The RVEs are defined by $[a_l, b_l]=[0,0.5^{L+1-l}]$ and the grid
sizes are $\mathfrak H=(0.25, 0.125, 0.0625)$. To ensure scale
separation even for the smallest RVE, we take 
$\dis \epsilon= \frac{b_1}{100}$. 

\begin{figure}[htp]\center
\subfigure[$\mathfrak M=(16M_3, 4M_3, M_3)$]{
\scalebox{.4}{ \includegraphics{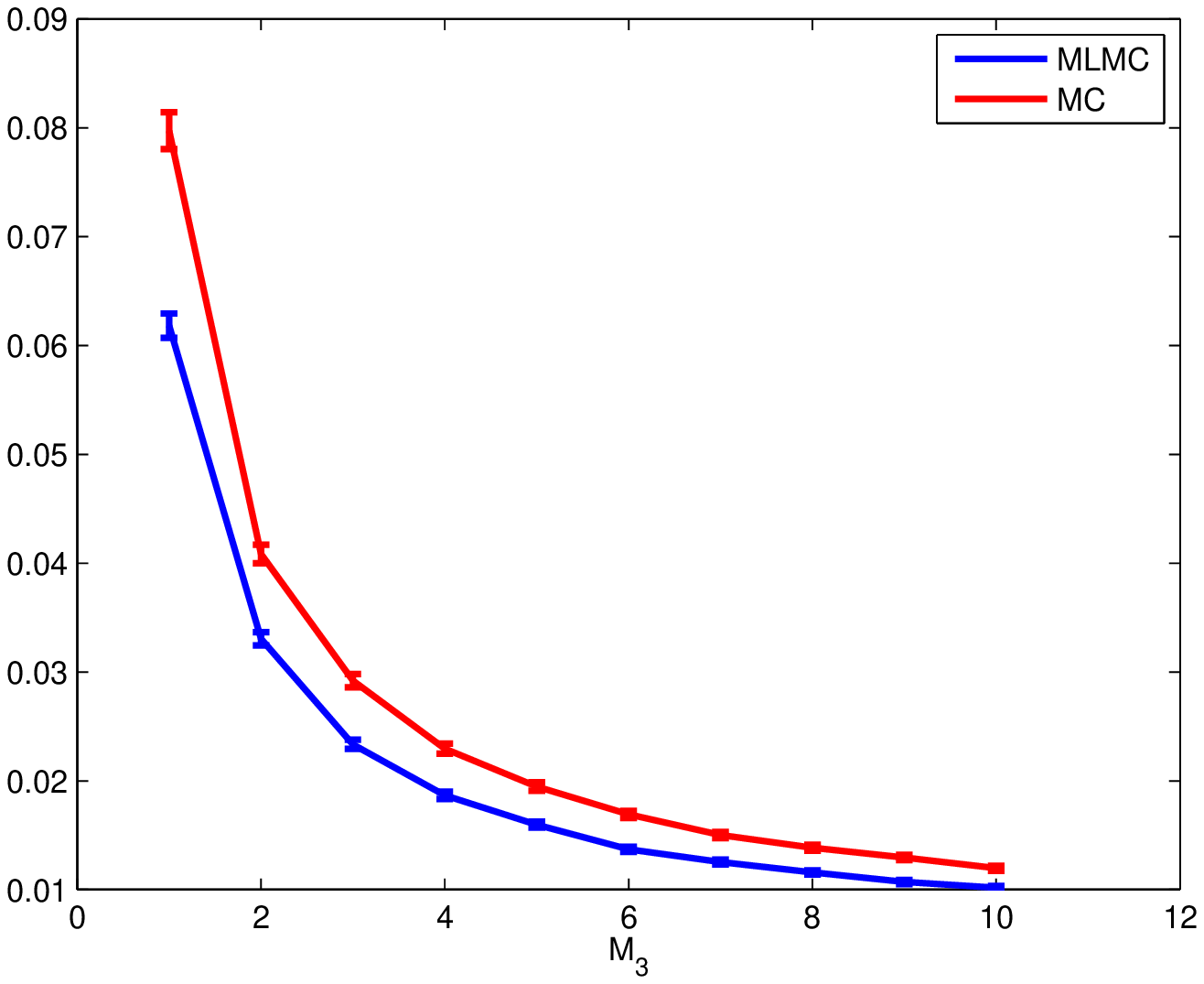}}}
\subfigure[$\mathfrak M=(4M_3, 2M_3, M_3)$]{
\scalebox{.4}{ \includegraphics{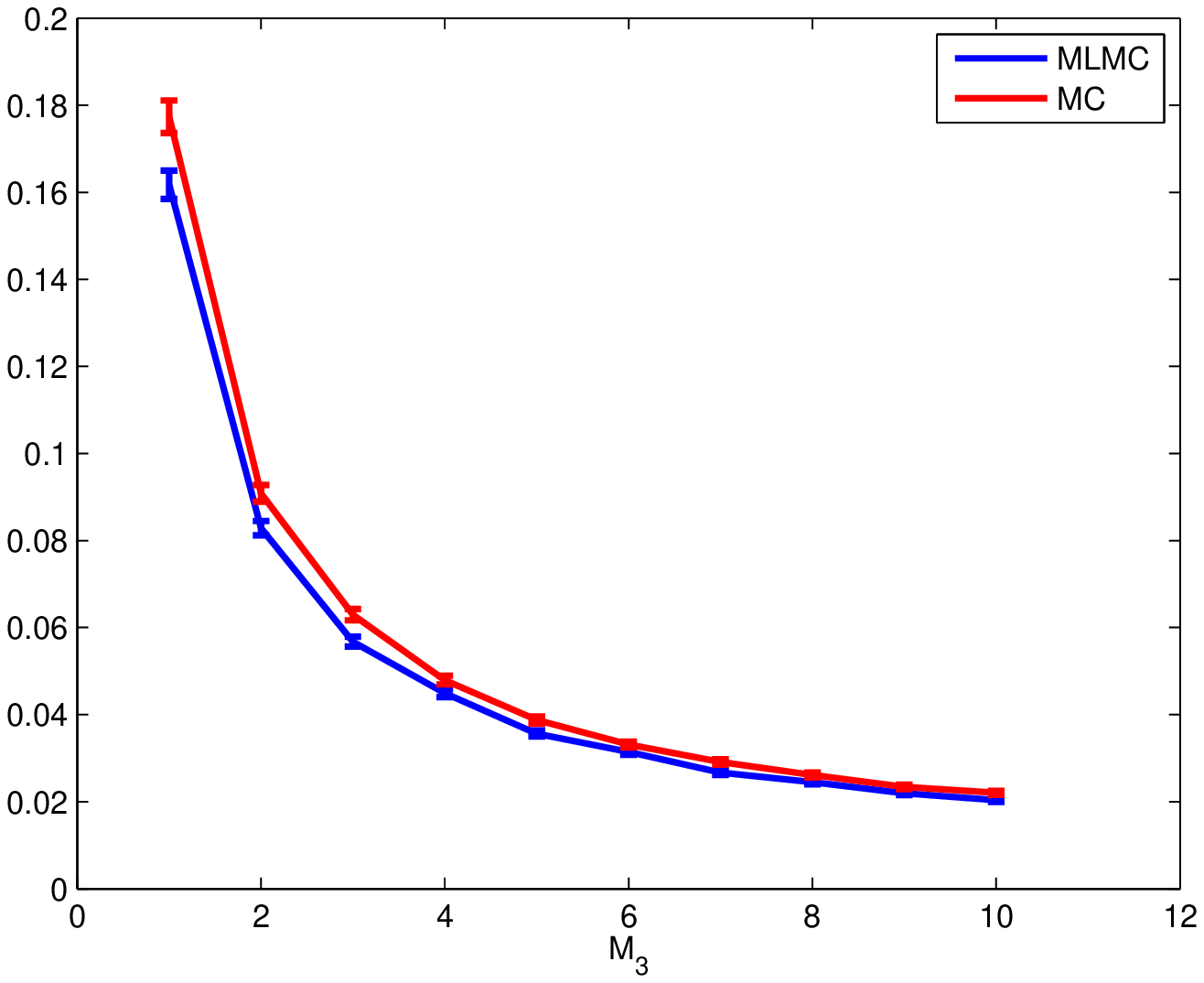}}}
\caption{Relative errors (in $L^2$ norm) of the solution (one
  dimensional example).}
\label{onedsolutioneps01}
\end{figure}

On Figure~\ref{onedsolutioneps01}, 
the accuracy of the MLMC approach is compared with that of the
MC approach at equal computational cost (error bars have been computed
using $Nb=20000$ different independent realizations of the whole
computation), for two choices of the number $\mathfrak M$ of
realizations at each level. We see that the choice $\mathfrak M=(16M_3,
4M_3, M_3)$, which is consistent with the rate $\beta=2$, yields the
best results (and an accuracy gain of 33 \%). 

\subsubsection{Two dimensional example}

We now turn to an example in dimension two. 
The reference problem~\eqref{microprob} is complemented with
homogeneous (zero) Dirichlet boundary conditions, and the 
source term is $f(x)=f(x_1,x_2)=100(x_1+x_2)$. 

\medskip

In the spirit of the Example 1 of Section~\ref{sec:coeff_2D}, we take
$$
A\left(x, \omega, \frac{x}{\epsilon}, \omega'\right) = 
\widetilde{A}(x,\omega) B\left(\frac{x}{\epsilon}, \omega'\right),
$$
where $\widetilde{A}$ and $B$ are scalar-valued, 
$B$ is a log-normal distributed random field, $B= e^K$, with $\EE(K) =
0$ and where the covariance function of $K(x,\omega')$ is 
$\text{cov}(x,x')= \sigma^2\exp \left(-\frac{| x-x'|^2}{\tau_0^2} \right)$, with
$\sigma = \tau_0 = \sqrt{2}$. The parameter $\epsilon$ is such that $\epsilon
\tau_0=0.04$. The macroscopic random field is given by
$$
\widetilde{A}(x,\omega)= 2 + 
\vert \omega_1 \sin(2\pi x_1)\vert
+\vert \omega_2 \sin(2\pi x_2)\vert 
+\vert \omega_3 \sin(\pi x_1)\vert
$$
with independent and normally distributed $\omega_k$, $1 \leq k \leq 3$. 

Since the coefficient is separable, we are in the setting described in
Section~\ref{sec:separable}. In particular, the RVE problems are
independent of the macroscopic point $x$, and we can use the MLMC
approach. For each level $1\leq l\leq L$, we hence solve the coarse
problem~\eqref{eq:homsoln} on a grid of size $H_l$, for $M_l$ realizations of
$\widetilde{A}$. This defines the solutions $u^k_l$, $1\leq k\leq M_l$,
$1 \leq l \leq L$. 

The MC approach consists in working only at the level $L$, and thus
solving, on a grid of size $H_L$, the problems
$$
-\text{div} \left( \widetilde{A}^k(x,\omega) E_{\widehat{m}}(B^*_L)
  \nabla u^k_L\right) = f \text{ in $D$}, \quad 1 \leq k \leq \widehat{M}.
$$
The reference solution is built as follows. At each level $l$, we first
solve~\eqref{eq:homsoln} with $m_l= \widetilde{m}$ and $1 \leq k \leq
M_l=\widetilde{M}$. The reference value is defined as the mean over both
the levels and the number of realizations of all these solutions:
$$
E^{\rm ref}_{\widetilde{M}, L} = \frac{1}{L} \sum_{l=1}^L
\frac{1}{\widetilde{M}} \sum_{k=1}^{\widetilde{M}} u_l^k.
$$
In practice, we take $\widetilde{M}= 1000$ and $\widetilde{m}= 50$.

We again work with $L=3$ different levels and we equate the
costs of the MC and the MLMC approaches for the computation of the
homogenized coefficients as well as that of the coarse scale
solutions. This respectively implies that the parameters of the MC
approach are 
$\dis \widehat{m} = \eta_3^{-2} 
\left(m_1 \eta_1^2+ m_2 \eta_2^2+ m_3 \eta_3^2 \right)$ 
and 
$\dis \widehat{M} = H_3^2 \left( M_3 H_3^{-2}+M_2 H_2^{-2}+ M_1 H_1^{-2}
  \right)$.

\begin{figure}[htp]
\center
\scalebox{.4}{ \includegraphics{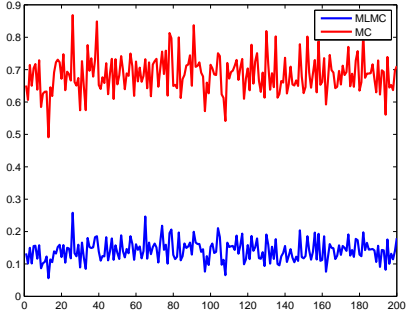}}
\caption{Relative $L_2$-errors $e_{MC}$ and $e_{MLMC}$ on the
  homogenized solution. We show the results for $200$ different independent
  realizations.} 
\label{solutionfig}
\end{figure}

On Figure~\ref{solutionfig}, we show the relative $L^2$-errors
\begin{eqnarray*}
e_{MLMC}(u_L) &=& 
\frac{\| E^{\rm ref}_{\widetilde{M},L}(u_L) - E^L(u_L) \|_{L^2(D)}}
{\| E^{\rm ref}_{\widetilde{M},L}(u_L) \|_{L^2(D)}}
\\ 
e_{MC}(u_L) &=& 
\frac{\| E^{\rm ref}_{\widetilde{M},L}(u_L) - E_{\widehat{M}}(u_L)
  \|_{L^2(D)}}
{\| E^{\rm ref}_{\widetilde{M},L}(u_L) \|_{L^2(D)}}
\end{eqnarray*}
computed with the parameters $\mathfrak M=(M_1,M_2,M_3)=(32, 32, 16)$
and $\mathfrak m=(m_1,m_2,m_3)=(50,40,20)$. 
Note that $\mathfrak M$ is chosen based on the calculations presented
in~\cite{BSZ10} (we have checked that these calculations also hold for
finite volume methods). 

We actually repeat the
whole procedure 200 times, and show on Figure~\ref{solutionfig} the 200
values of the relative errors that we found. We see that these errors
are essentially the same for all the realizations. A gain in accuracy of
the order of 5 is obtained when using the MLMC approach, for an equal
cost:
$$
\EE(e_{MLMC}) \approx 0.1411, 
\quad
\EE(e_{MC}) \approx 0.6851.
$$
The standard deviation of the MLMC error is also smaller:
$$
{\rm std}_{MLMC} = 0.0324,
\quad
{\rm std}_{MC} = 0.0565.
$$

\section*{Acknowledgments}

A part of this work was done while YE was visiting ENPC and ITWM. YE is
grateful for the support from ENPC and INRIA as well as Humboldt
Foundation and ITWM.
The research of CK was partially supported by the DFG Project IL
55/1-2.
The work of FL is partially supported by ONR under Grant
N00014-12-1-0383. 
FL warmly thanks the Fraunhofer Institute ITWM (where this work was
initiated) and the Texas A~\&~M  University for their kind
hospitalities. 

\appendix

\section{Appendix: weighted MLMC approach analysis}
\label{sec:app}

We estimate here the error associated to the weighted MLMC approach
introduced in Section~\ref{solutionGeneral}. To this aim, it
is useful to introduce the function
$$
\widetilde{u} = \sum_{l=1}^L \frac{\alpha_l}{M_l} \sum_{j=l}^L
(M_j-M_{j+1}) \left(u^*_{\eta_j,H_l}-u^*_{\eta_j,H_{l-1}}\right).
$$
We indeed note that 
$\EE\left( E^{L*}_{weighted} \right) = \EE(\widetilde{u})$. The error
between the computed quantity $E^{L*}_{weighted}$ and the exact quantity
$\EE(u^*)$ is thus composed of a statistical error (the expectation
of $E^{L*}_{weighted}$ is only approximately estimated) and of a systematic
error, due to the fact that 
$\EE\left( E^{L*}_{weighted} \right) = \EE(\widetilde{u}) \neq
\EE(u^*)$. We successively estimate these two contributions. 

\paragraph{Systematic error estimation}

Following the same lines as in Section~\ref{sec:mlmc}, we obtain that 
$$
\left\| \alpha_1 u^* -\widetilde{u} \right\|
\leq
\sum_{l=1}^L {\cal E}_{l,l} \left[ \alpha_l -\alpha_{l+1} \right]
$$
where we have set $\alpha_{L+1} = 0$ and 
$\dis {\cal E}_{j,l} := \left\| u^*- u^*_{\eta_j,H_l}
\right\|$. Choosing now 
\begin{equation}
\label{eq:choice_w-mlmc}
\alpha_l = \sum_{j=l}^L \widetilde{\alpha}_j 
\frac{{\cal E}_{L,L}}{{\cal E}_{j,j}}
\end{equation}
to equilibrate the terms in the above error bound, we get
$$
\left\| \alpha_1 u^* -\widetilde{u} \right\|
\leq
{\cal E}_{L,L} \sum_{l=1}^L \widetilde{\alpha}_l
\quad \text{with} \quad
\alpha_1 = \sum_{j=1}^L \widetilde{\alpha}_j 
\frac{{\cal E}_{L,L}}{{\cal E}_{j,j}}.
$$
As shown in Section~\ref{sec:separable}, we have
$\dis {\cal E}_{j,l} \lesssim H_l + 
\left( \frac{\epsilon}{\eta_j} \right)^{\beta/2} + \frac{1}{\sqrt{m_j}}$. 

\medskip

For the standard MC approach, the systematic error reads
$$
\| u^* - u_{\widehat{\eta}, \widehat{H}} \| \lesssim
\widehat{H} + \widehat{\delta}.
$$
To have the same systematic error, we choose the coarse grid size
$\dis \widehat{H} = \sum_{l=1}^L (\alpha_l-\alpha_{l+1}) H_l$ and RVEs
of size $\widehat {\eta}$ so that
$\dis \widehat{\delta} 
= 
\left( \frac{\epsilon}{\widehat{\eta}} \right)^{\beta/2}
=
\sum_{l=1}^L (\alpha_l-\alpha_{l+1}) \delta_l$.

\paragraph{Statistical error estimation}

The statistical error of the weighted MLMC approximation satisfies
$$
\left\| \EE(\widetilde{u}) - 
\sum_{l=1}^L \alpha_l E^*_{M_l}(u^*_l -u^*_{l-1}) \right\|
\lesssim
\sum_{l=1}^L \frac{\alpha_l}{\sqrt{M_l}} 
(H_l+\delta_l) + \frac{\alpha_1}{\sqrt{M_1}}.
$$
To equate the error terms in the above sum, we choose
\begin{equation}
\label{eq:choix_Ml}
M_l= C \left(\frac{\alpha_l(H_l+\delta_l)}
{\gamma_l(\widehat{H}+ \widehat{\delta})}\right)^2 
\ \text{for $l\geq2$},
\quad
M_1= C \left(\frac{\alpha_1}
{\gamma_1(\widehat{H}+ \widehat{\delta})}\right)^2,
\end{equation}
for some constant $C$ and some parameters $\gamma_l$.
The statistical error then satisfies
$$
\left\| \EE\left(E^{L*}_{weighted}\right) - E^{L*}_{weighted} \right\|
= 
O(\widehat{H} + \widehat{\delta}).
$$
For the MC approach, we choose $\widehat{M}= C
(\widehat{H}+\widehat{\delta})^{-2}$ independent realizations and thus
get a statistical error of the same order. 

\paragraph{Cost comparison}

Now that we have chosen parameters such that the MC and the weighted
MLMC approaches share the same accuracy, we are in position to compare
their cost. 

\begin{figure}
 [htp]\center
\subfigure[$\epsilon=\frac{\eta_1}{10}$]{
\scalebox{.4}{ \includegraphics{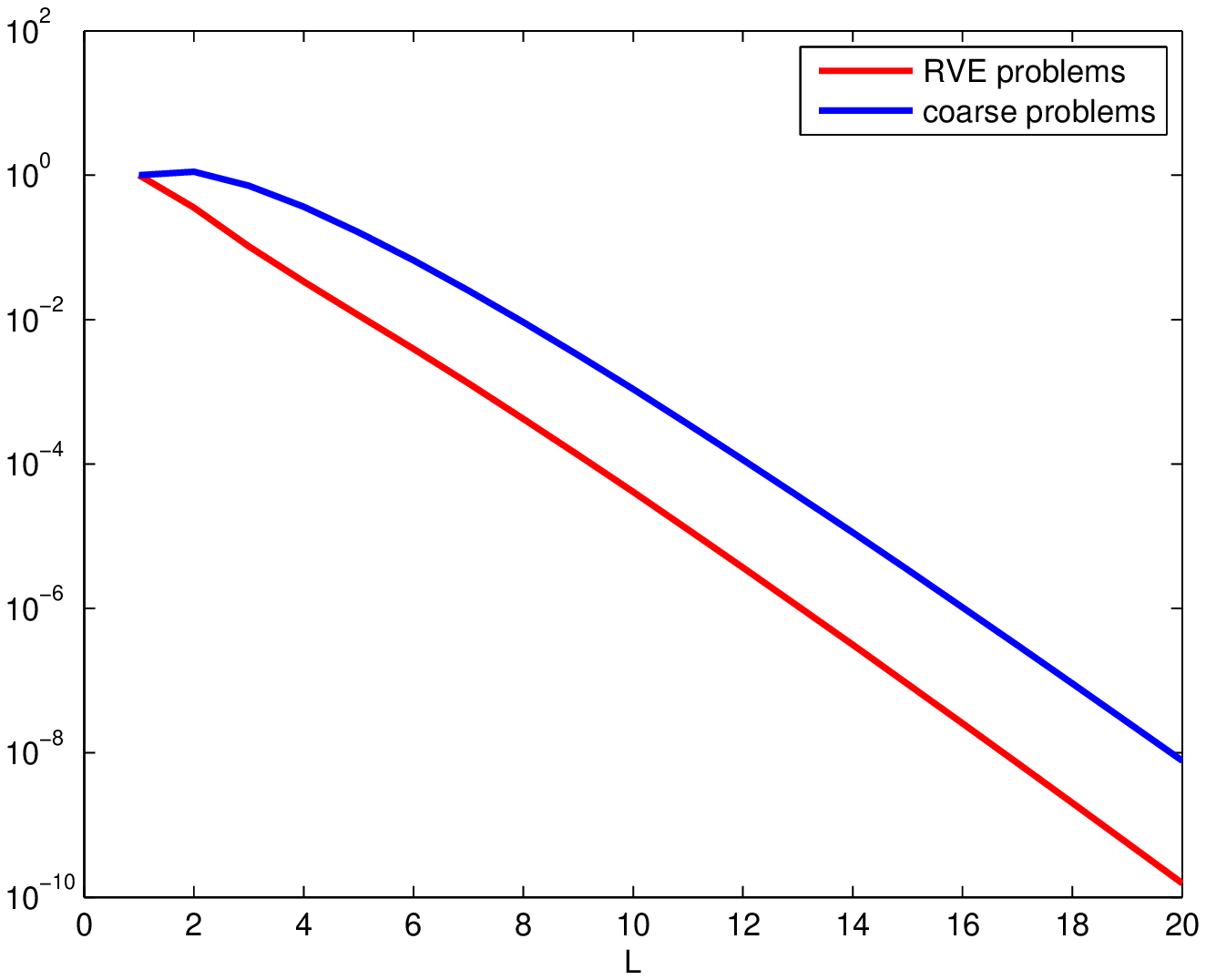}}}
\subfigure[$\epsilon=\frac{2^{-50}}{10}$]{
\scalebox{.4}{ \includegraphics{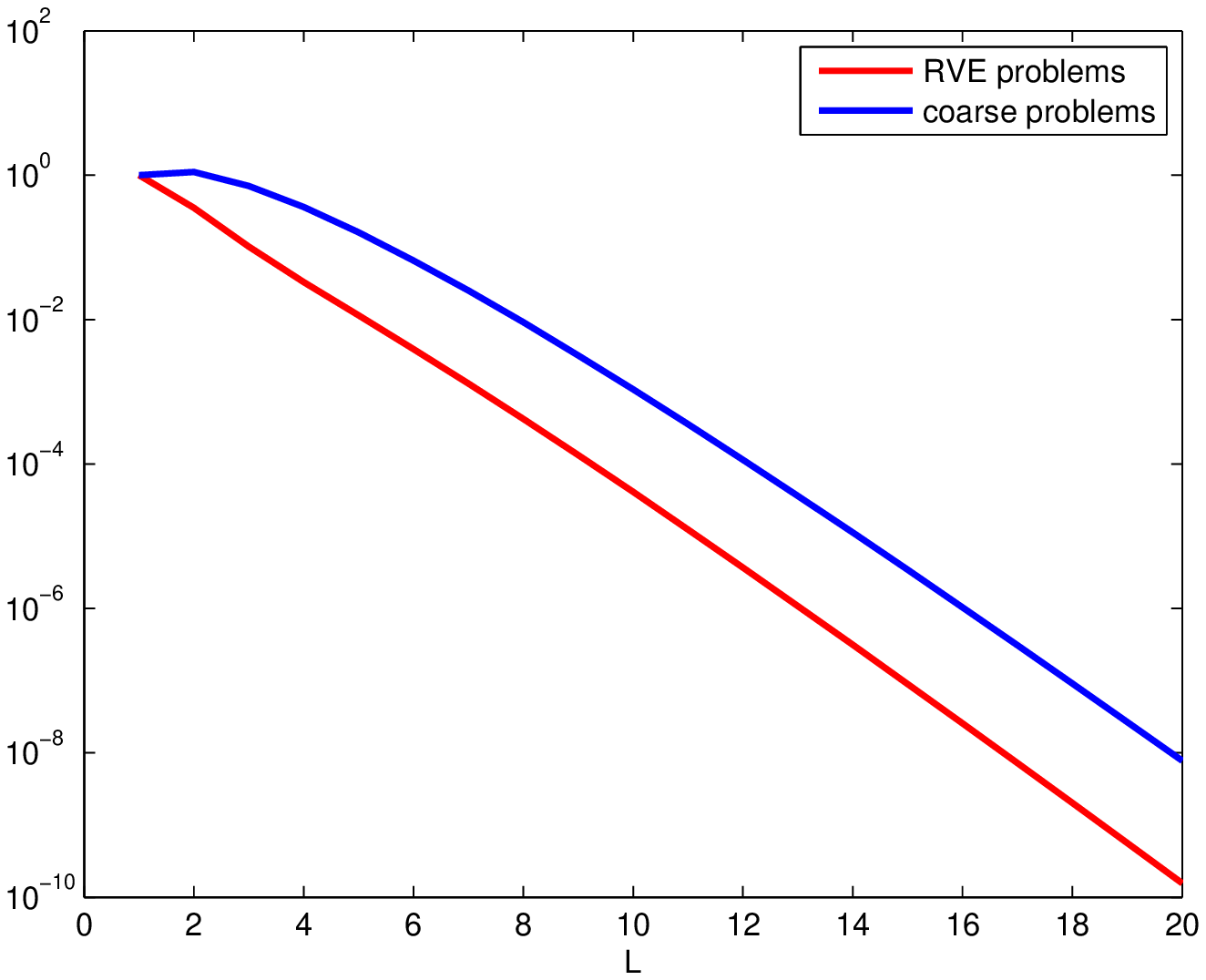}}}
\caption{Work ratios $\frac{W^{w-MLMC}_{\rm RVE}}{W^{MC}_{\rm RVE}}$ and
  $\frac{W^{w-MLMC}_{\rm coarse}}{W^{MC}_{\rm coarse}}$ as a function of
  $L$. We work with $\beta=2$, $\eta_l=2^{l-L}$ and $\gamma_l = 1/L$. We
  choose $M_l$ according to~\eqref{eq:choix_Ml} and $\alpha_l$ according
  to~\eqref{eq:choice_w-mlmc} with 
  $\widetilde{\alpha}_j = C$, where the constant $C$ is such that
  $\alpha_1=1$.}
\label{wmlmcworkratio}
\end{figure}

As above, the cost of solving the coarse scale problems is
$$
W^{w-MLMC}_{\rm coarse} = \sum_{l=1}^L M_l H_l^{-2}
\quad \text{and} \quad
W^{MC}_{\rm coarse} = \widehat{M} \widehat{H}^{-2}.
$$
The dominating part of the computational cost however lies in 
solving the local RVE problems. For the MC approach, we assume that we
need to solve these problems at $\widehat{N}$ macroscopic points $x$. We
thus have
$$
W^{MC}_{\rm RVE} =\widehat{M} \widehat{N} 
\left(\frac{\widehat{\eta}}{\epsilon}\right)^2
= C \frac{\widehat{N} \widehat{\eta}^2}
{\epsilon^2 (\widehat{H} + \widehat{\delta})^2}.
$$
For the weighted MLMC approach, we assume that, at each level $l$, 
we solve local RVE problems at $N_l \leq H_l^{-2}$ macroscopic points
(with $N_1<N_2<\cdots<N_L$). At each of these points, we only need to
consider $M_l-M_{l+1}$ realizations. The computational work thus reads
$$
W^{w-MLMC}_{\rm RVE}
= 
\sum_{l=1}^L (M_l-M_{l+1})\left(\frac{\eta_l}{\epsilon}\right)^2 N_l.
$$
On Figure~\ref{wmlmcworkratio}, we show the ratio of the works for
solving the coarse problems and the RVE problems, as a function of the
number of levels $L$. The figure is made with the parameter $\beta=2$
(which corresponds to a Central Limit Theorem type convergence, see discussion
below~\eqref{eq:hyp}). As on Figure~\ref{rvework}, we consider two
possible regimes for $\epsilon$. We see that a significant gain is
achieved even for moderate values of $L$. 

\end{document}